\documentclass[12pt]{article}
\usepackage{latexsym, epsfig, amssymb, amsmath}
\usepackage[dvips]{color}
\usepackage{graphicx}
\usepackage{anysize}

\makeatletter

\renewcommand{\baselinestretch}{1.3}
\marginsize{1 in}{1 in}{1 in}{1 in}

\def\singlespace{\def\baselinestretch{1}\@normalsize}

\newtheorem{theorem}{Theorem}

\newtheorem{lemma}{Lemma}

\def\Q{{\mathbb Q}}        
\def\Z{{\mathbb Z}}        
\def\P{{\mathbb P}}        
\def\1{{\mathbf 1}}        

\def\a{{\boldsymbol{\alpha}}}
\def\b{{\boldsymbol{\beta}_n}}
\def\bhat{{\hat{\boldsymbol{\beta}}_n}}
\def\v{{\mathbf{v}}}
\def\X{{\mathbf{X}}}
\def\Z{{\mathbf{Z}}}
\def\Q{{\hat{Q}_n}}
\def\g{{\gamma_n}}
\def\gsq{{\gamma_n^2}}
\def\gcu{{\gamma_n^3}}
\def\ahat{{\hat{\boldsymbol{\alpha}}_{\boldsymbol{\beta}_n}}}
\def\ahato{{\hat{\boldsymbol{\alpha}}_{\boldsymbol{\beta}_{n0}}}}
\def\ahatp{{\hat{\boldsymbol{\alpha}}_{\boldsymbol{\beta}_n}^{\prime}}}
\def\ahatop{{\hat{\boldsymbol{\alpha}}_{\boldsymbol{\beta}_{n0}}^{\prime}}}

\def\ao{{\boldsymbol{\alpha}_{\boldsymbol{\beta}_{n0}}}}
\def\ap{{\boldsymbol{\alpha}_{\boldsymbol{\beta}_n}^{\prime}}}
\def\aop{{\boldsymbol{\alpha}_{\boldsymbol{\beta}_{n0}}^{\prime}}}

\def\abar{{\bar{\boldsymbol{\alpha}}_{ni}}}
\def\bstar{{\boldsymbol{\beta}^*}}
\def\bhatstar{{\hat{\boldsymbol{\beta}}^*}}
\def\gvec{{\boldsymbol{\gamma}_n}}
\def\ghatvec{{\hat{\boldsymbol{\gamma}}_n}}
\def \bbeta{\boldsymbol{\beta}}
\def \lh{\log^{1/2}(1/h)}
\def \bB{\mbox{\bf B}}
\def \bX{\mbox{\bf X}}
\def \bZ{\mbox{\bf Z}}

\def\pcv{\stackrel{\scriptscriptstyle \P}{\longrightarrow}}        
\def\dcv{\stackrel{\scriptscriptstyle \mathcal{D}}{\longrightarrow}}        

\makeatletter
\def\table{\@ifnextchar[{\table@i}{\table@i[\fps@table]}}
\def\table@i[#1]{\@float{table}[#1]\footnotesize}
\makeatother

\begin{document}

\title{\vspace{-0.9 in}
\bf Profile-Kernel Likelihood Inference With Diverging Number of
Parameters
\thanks{Clifford Lam is PhD student, Department of Operation Research and
Financial Engineering, Princeton University, Princeton, NJ 08544
(email: wlam@princeton.edu); Jianqing Fan is Professor, Department
of Operation Research and Financial Engineering, Princeton
University, Princeton, NJ 08544 (email: jqfan@princeton.edu).
Financial support from the NSF grant DMS-0354223, DMS-0704337 and
NIH grant R01-GM072611 is gratefully acknowledged.}
\date{}
\author{By Clifford Lam and Jianqing Fan    \\
Department of Operations Research and Financial Engineering \\
Princeton University, Princeton, NJ, 08544}} \maketitle

\vspace{-0.25 in}
\begin{singlespace}
\begin{quotation}
\indent

The generalized varying coefficient partially linear model with
growing number of predictors arises in many contemporary
scientific endeavor.  In this paper we set foot on both
theoretical and practical sides of profile likelihood estimation
and inference. When the number of parameters grows with sample
size, the existence and asymptotic normality of the profile
likelihood estimator are established under some regularity
conditions. Profile likelihood ratio inference for the growing
number of parameters is proposed and Wilk's phenomenon is
demonstrated. A new algorithm, called the accelerated
profile-kernel algorithm, for computing profile-kernel estimator
is proposed and investigated. Simulation studies show that the
resulting estimates are as efficient as the fully iterative
profile-kernel estimates. For moderate sample sizes, our proposed
procedure saves much computational time over the fully iterative
profile-kernel one and gives stabler estimates. A set of real data
is analyzed using our proposed algorithm.
\end{quotation}
\end{singlespace}

{\em Short Title}: High-dimensional profile likelihood.

{\em AMS 2000 subject classifications}. Primary 62G08; secondary
62J12, 62F12.

{\em Key words and phrases}. Generalized linear models, varying
coefficients, high dimensionality, asymptotic normality, profile
likelihood, generalized likelihood ratio tests.

\newpage

\section{Introduction}

Semiparametric models with large number of predictors arise
frequently in many contemporary statistical studies.  Large data
set and high-dimensionality characterize many contemporary
scientific endeavors (\cite{Donoho(2000)}; \cite{FanLi(2006)}).
Statistical models with many predictors are frequently employed to
enhance the explanatory and predictive powers.  At the same time,
semiparametric modeling is frequently incorporated to balance
between modeling biases and ``curse of dimensionality''. Profile
likelihood techniques (\cite{MurphyvanderVaart(2000)}) are
frequently applied to this kind of semiparametric models.  When
the number of predictors is large, it is more realistic to regard
it growing with the sample size. Yet, few results are available
for semiparametric profile inferences when the number of
parameters diverges with sample size.  This paper focuses on
profile likelihood inferences with diverging number of parameters
in the context of the generalized varying coefficient partially
linear model (GVCPLM).

GVCPLM is an extension the generalized linear model
(\cite{McCullaghNelder(1989)}) and the generalized
varying-coefficient model (\cite{HastieTibshirani(1993)};
\cite{Caietal(2000)}).  It allows some coefficient functions to
vary with certain covariates $U$ such as age
(\cite{FanPeng(2004)}), toxic exposure level or time variable in a
longitudinal data or survival analysis (\cite{Murphy(1993)}).
Therefore, general interactions, not just the linear interaction
as in parametric models, between the variable $U$ and these
covariates are explored nonparametrically.

If Y is a response variable and $(U, \mathbf{X},\mathbf{Z})$ is
the associated covariates, then by letting
$\mu(u,\mathbf{x},\mathbf{z}) =
E\{Y|(U,\mathbf{X},\mathbf{Z})=(u,\mathbf{x}, \mathbf{z})\}$, the
GVCPLM takes the form
\begin{equation}\label{eqn:GVCPLM}
g\{\mu(u,\mathbf{x},\mathbf{z})\} = \mathbf{x}^T\a(u) +
\mathbf{z}^T \boldsymbol{\beta},
\end{equation}
where $g(\cdot)$ is a known link function, $\bbeta$ a vector of
unknown regression coefficients and $\a(\cdot)$ a vector of
unknown regression functions. One of the advantages over the
varying coefficient model is that GVCPLM allows more efficient
estimation when some coefficient functions are not really varying
with $U$, after adjustment of other genuine varying effects. It
also allows more interpretable model, where primary interest is
focused on the parametric component.

\subsection{A motivating example}\label{subsect:example}
We use a real data example to demonstrate the need for GVCPLM. The
Fifth National Bank of Springfield faced a gender discrimination
suit in which female received substantially smaller salaries than
male employees. This example is based on a real case with data
dated 1995. Only the bank's name is changed.  See Example 11.3 of
\cite{Albrightetal(1999)}. Among 208 employees, eight variables
are collected.  They include employee's salary; age; year hired;
number of years of working experience at another bank; gender; PC
Job, a dummy variable with value 1 if the employee's job is
computer related; educational level, a categorical variable with
categories 1 (finished school), 2 (finished some college courses),
3 (obtained a bachelor's degree), 4 (took some graduate courses),
5 (obtained a graduate degree); job grade, a categorical variable
indicating the current job level, the possible levels being 1--6
(6 the highest).

\cite{FanPeng(2004)} has conducted such a salary analysis using an
additive model with quadratic spline and does not find a
significant evidence of gender difference.  However, salary is
directly related to the job grade.  With the adjustment for the
job grade, the salary discrimination can not easily be seen.  An
important question then arises if female employees have lower
probability getting promoted. In analyzing such probability, a
common tool will be the logistic regression, a class of the
generalized linear model (for example, see
\cite{McCullaghNelder(1989)}).

To this end, we create a binary response variable
\textbf{HighGrade4}, indicating if Job Grade is greater than 4.
The associated covariates are \textbf{Female}(1 for female
employee and 0 otherwise), \textbf{Age},
\textbf{TotalYrsExp}(total years of working experience),
\textbf{PCJob}, \textbf{Edu}(level of education). If the covariate
\textbf{Female} has a significantly negative coefficient, then it
would suggest that female employees are harder to promote to
higher grade jobs.

However, in a simple logistic regression, the effect of a
covariate cannot change with another covariate nonparametrically.
Table \ref{table:example}  shows the proportion of employees
having a job grade higher than 4, categorized by \textbf{Age} and
\textbf{TotalYrsExp}. Clearly interactions between \textbf{Age}
and \textbf{TotalYrsExp} have to be considered.

\begin{table}
  \centering
  \caption{Proportions of employees having job grade higher than 4}
  \begin{tabular}{cccc}\label{table:example}\\
  \hline
       &   \multicolumn{3}{c}{Covariate \textbf{TotalYrsExp}} \\
       &   0-7     &    8-16     & $\geq$17        \\
  \hline
 Age $\leq$ 35 & 1/11 & 1/9 & 0        \\
 Age $>$ 35      & 0 & 2/11 & 8/21      \\
  \hline
\end{tabular}
\end{table}

This can be done by creating categorical variables over the
covariate \textbf{Age}. However this would increase the number of
predictors considerably if we create many categories of
\textbf{Age}. More importantly, we do not know where to draw the
borders of each \textbf{Age} category and how many categories
should be produced. This problem is nicely overcome if we allow
the coefficient of \textbf{TotalYrsExp} to vary with \textbf{Age},
so that we obtain a coefficient function of \textbf{Age} for
\textbf{TotalYrsExp}. See section \ref{subsect:realdataanalysis}
for a detail analysis of the data.

If interactions between different variables are considered, then
the number of predictors will be large compare with the sample
size $n=208$. This motivates us to consider the setting $p_n
\rightarrow \infty$ as $n \rightarrow \infty$ and present general
theories in section \ref{sect:properties}, where such a setting
will be faced by many modern statistical applications.

\subsection{Goals of the paper}

When the number of parameters $\bbeta$ is fixed and the link $g$
is identity, the model (\ref{eqn:GVCPLM}) has been considered by
\cite{Zhangetal(2002)}, \cite{Lietal(2002)} and
\cite{Xiaetal(2004)}, and \cite{Ahmadetal(2005)}.
\cite{FanHuang(2005)} propose a profile-kernel inference for such
a varying coefficient partial linear model (VCPLM) and
\cite{LiLiang(2005)} considered a backfitting-based procedure for
model selection in VCPLM. All of these papers rely critically on
the explicit form of the estimation procedures and the techniques
can not easily be applied to the GVCPLM.

Modern statistical applications often involve estimation of large
number of parameters. It is of interest to derive asymptotic
properties for the profile likelihood estimators under model
(\ref{eqn:GVCPLM}) when number of parameters diverges.  The
fundamental questions arise naturally whether the profile
likelihood estimator (\cite{MurphyvanderVaart(2000)}) still
possesses efficient sampling properties, whether the profile
likelihood ratio test for the parametric component possesses Wilks
type of phenomenon, namely whether the asymptotic null
distributions are independent of nuisance functions and
parameters, and whether the usual sandwich formula provides a
consistent estimator of the covariance matrix of the profile
likelihood estimator. These questions are poorly understood and
will be thoroughly investigated in Section 2. Pioneering work on
statistical inference with diverging number of parameters include
\cite{Huber(1973)} which gave related results on M-estimators, and
\cite{Portnoy(1988)} which analyzed a regular exponential family
under the same setting. \cite{FanPeng(2004)} studied the penalized
likelihood approach under such setting, whereas
\cite{Fanetal(2005)} investigated a semiparametric model with
growing number of nuisance parameters.

Another goal of this paper is to provide an efficient algorithm
for computing profile likelihood estimates under the model
(\ref{eqn:GVCPLM}). To this end, we propose a new algorithm,
called the accelerated profile-kernel algorithm, based on an
important modification of the Newton-Raphson iterations.
Computational difficulties (\cite{LinCarroll(2006)}) of the
profile-kernel approach is significantly reduced, while nice
sampling properties of such approach over the backfitting
algorithm (e.g. \cite{Huetal(2004)}) are retained. This will be
convincingly demonstrated in Section \ref{sect:simulation}, where
the Poisson and Logistic specifications are considered for
simulations. A new difference-based estimate for the parametric
component is proposed as an initial estimate of our proposed
profile-kernel procedure. Our method expands significantly the
idea used in \cite{Yatchew(1997)} and \cite{FanHuang(2005)} for
the partial linear model.

 The outline of the paper is as follows. In Section
\ref{sect:properties} we briefly introduce the profile likelihood
estimation with local polynomial modeling and present our main
asymptotic results. Section \ref{sect:computation} turns to the
computational aspect,  discussing the elements of computing in the
accelerated profile-kernel algorithm.  Simulation studies and an
analysis of real data set are given Section \ref{sect:simulation}.
The proofs of our results are given in Section \ref{sect:proof},
and technical details in the appendix.

\section{Properties of profile likelihood
inference}\label{sect:properties} \setcounter{equation}{0}

 Let $(Y_{ni} ; \textbf{X}_i, \textbf{Z}_{ni}, U_i)$, where
$1 \leq i \leq n$ be a random sample where $Y_{ni}$ is a scalar
response variable, $U_i$, $\textbf{X}_i \in \mathbb{R}^q$ and $
\textbf{Z}_{ni} \in \mathbb{R}^{p_n}$ are vectors of explanatory
variables. We consider model (\ref{eqn:GVCPLM}) with $\b$ and
$\mathbf{Z}_n$ having dimensions $p_n \rightarrow \infty$ as $n
\rightarrow \infty$. Like the distributions in the exponential
family, we assume that the conditional variance depends on the
conditional mean so that $\text{Var}(Y|U, \bX, \bZ_n) =
V(\mu(u,\bX, \bZ_n))$ for a given function $V$ (Our result is
applicable even when $V$ is multiplied by an unknown scale). Then,
the conditional quasi-likelihood function is given by
$$ Q(\mu, y) = \int_{\mu}^y \frac{s-y}{V(s)}ds. $$
As in \cite{SeveriniWong(1992)}, we denote by $\a_\b(u)$ the
`least favorable curve' of the nonparametric function $\a(u)$,
which is defined as the one that maximizes
\begin{equation}\label{eqn:defining1}
E_0\{Q(g^{-1}(\boldsymbol{\eta}^T\X + \b^T\Z_n), Y_n)|U=u \}
\end{equation}
with respect to $\boldsymbol{\eta}$,  where $E_0$ is the
expectation taken under the true parameters $\a_0(u)$ and $\b_0$.
As will be discussed in section \ref{subsect:asympnorm}, through
the use of least favorable curve, no undersmoothing of the
nonparametric component is required to achieve asymptotic
normality when $p_n$ is diverging with $n$. Note that
$\a_{\b_0}(u) = \a_0(u) $. Under some mild conditions, it
satisfies
\begin{equation}\label{eqn:defining}
\frac{\partial}{\partial\boldsymbol{\eta}}
E_0\{Q(g^{-1}(\boldsymbol{\eta}^T\X + \b^T\Z_n), Y_n)|U=u
\}|_{\boldsymbol{\eta} = \a_\b(u)}=0.
\end{equation}
The profile-likelihood function for $\b$ is then
\begin{equation}\label{eqn:glf}
Q_n(\b) = \sum_{i=1}^{n}Q\{g^{-1}(\a_{\b}(U_i)^{T}\mathbf{X}_i +
\mathbf{\boldsymbol{\beta}}_n^{T}\mathbf{Z}_{ni}), Y_{ni}\},
\end{equation}
if the least-favorable curve $\a_\b(\cdot)$ is known.

The least-favorable curve defined by (\ref{eqn:defining1}) can be
estimated by its sample version through a local polynomial
regression approximation.  For $U$ in a neighborhood of $u$,
approximate the $j^{th}$ component of $\a_\b(\cdot)$ as
\begin{eqnarray*}
  \alpha_{j}(U) &\approx& \alpha_{j}(u) + \frac{\partial\alpha_{j}(u)}{\partial{u}}(U-u) + \cdots + \frac{\partial^p\alpha_{j}(u)}{\partial{u^p}}(U-u)^p/p!\\
                &\equiv& a_{0j} + a_{1j}(U-u) + \cdots +
                a_{pj}(U-u)^p/p!.
\end{eqnarray*}
Denoting $\mathbf{a_r} = (a_{r1},\cdots, a_{rq})^T $ for
$r=0,\ldots,p$, for each given $\b$, we then maximize the local
likelihood
\begin{equation}\label{eqn:llf}
\sum_{i=1}^n
Q\{g^{-1}(\sum_{r=0}^{p}\mathbf{a_r}^T\mathbf{X}_i(U_i - u)^r/r! +
\boldsymbol{\beta}_n^{T}\mathbf{Z}_{ni}), Y_{ni}\}K_h(U_i - u)
\end{equation}
with respect to $\mathbf{a_0},\cdots, \mathbf{a_p}$, where
$K(\cdot)$ is a kernel function and $K_h(t) = K(t/h)/h$ is a
re-scaling of $K$ with bandwidth $h$. Thus, we get estimate
$\hat{\mathbf{\boldsymbol{\alpha}}}_{\boldsymbol{\beta}_n}(u) =
\mathbf{\hat{a}_0}(u)$.

Plugging our estimates into the profile-kernel likelihood function
(\ref{eqn:glf}), we have
\begin{equation}\label{eqn:glfhat}
\hat{Q}_n(\boldsymbol{\beta}_n) = \sum_{i=1}^n
Q\{g^{-1}(\hat{\boldsymbol{\alpha}}_{\boldsymbol{\beta}_n}(U_i)^{T}\mathbf{X}_i
+ \boldsymbol{\beta}_n^{T}\mathbf{Z}_{ni}), Y_{ni}\}.
\end{equation}
Maximizing $\hat{Q}_n(\boldsymbol{\beta}_n)$ with respect to
$\boldsymbol{\beta}_n$ to get $\hat{\boldsymbol{\beta}}_n$. With
$\hat{\boldsymbol{\beta}}_n$, the varying coefficient functions
are estimated as
$\hat{\boldsymbol{\alpha}}_{\hat{\boldsymbol{\beta}}_n}(u)$.

One property of the profile quasi-likelihood is that the first and
second order Bartlett's identities continue to hold. In
particular, with the definition given by (\ref{eqn:glf}), then for
any $\b$, we have
\begin{equation}\label{eqn:bartlett}
\mathbf{E}_\b\biggl(\frac{\partial{Q_n}}{\partial{\boldsymbol{\beta}_n}}\biggr)=0,
\qquad
\mathbf{E}_\b\biggl(\frac{\partial{Q_n}}{\partial{\boldsymbol{\beta}_n}}\frac{\partial{Q_n}}{\partial{\boldsymbol{\beta}_n^T}}\biggr)
= -\mathbf{E}_\b\biggl(
\frac{\partial^2Q_n}{\partial\boldsymbol{\beta}_n\partial{\boldsymbol{\beta}_n^T}}
\biggr).
\end{equation}
See \cite{SeveriniWong(1992)} for more details. These properties
give rise to the asymptotic efficiency of the profile likelihood
estimator.

\vspace{12pt}
\subsection{Consistency and asymptotic normality of $\bhat$}\label{subsect:asympnorm}

We need Regularity Conditions (A) - (G) in Section
\ref{sect:proof} for the following results.

\begin{theorem}\label{thm:ple} \emph{(Existence of profile
likelihood estimator).} Assume that Conditions (A)-(G) are
satisfied. If $p_n^4/n \rightarrow 0$ as $n \rightarrow \infty$
and $h=O(n^{-a})$ with $(4(p+1))^{-1} < a < 1/2$, then there is a
local maximizer $\hat{\boldsymbol{\beta}}_n \in \Omega_n$ of
$\hat{Q}_n(\boldsymbol{\beta}_n)$ such that
$\|\hat{\boldsymbol{\beta}}_n - \boldsymbol{\beta}_{n0}\| =
O_P(\sqrt{p_n/n})$.
\end{theorem}

The above rate is the same as the one established by
\cite{Huber(1973)} for the M-estimator.

Note that the optimal bandwidth $h=O(n^{-1/{(2p+3)}})$ is included
in Theorem \ref{thm:ple}. Hence $\sqrt{n/p_n}$-consistency is
achieved without the need of undersmoothing of the nonparametric
component. In particular, when $p_n$ is fixed, the result is in
line with those, for instance, by \cite{SeveriniStaniswalis(1994)}
in a different context.

Define $I_n(\b) = n^{-1}
\mathbf{E}_\b(\frac{\partial{Q_n}}{\partial{\boldsymbol{\beta}_n}}\frac{\partial{Q_n}}{\partial{\boldsymbol{\beta}_n^T}})
$, which is an extension of the Fisher matrix. Since the
dimensionality grows with sample size, we need to consider the
arbitrary linear combination of the profile kernel estimator
$\hat{\boldsymbol{\beta}}_n$ as stated in the following theorem.

\begin{theorem}\label{thm:asymnorm} \emph{(Asymptotic
normality).} Under Conditions (A) - (G), if $p_n^5/n = o(1)$ and
$h=O(n^{-a})$ for $3/(10(p+1)) < a < 2/5$, then the consistent
estimator $\hat{\boldsymbol{\beta}}_n$ in Theorem \ref{thm:ple}
satisfies
$$
\sqrt{n}A_nI_n^{1/2}(\boldsymbol{\beta}_{n0})(\hat{\boldsymbol{\beta}}_n
- \boldsymbol{\beta}_{n0}) \dcv N(0, G), $$ where $A_n$ is an $l
\times p_n$ matrix such that $A_nA_n^T \rightarrow G$, and $G$ is
an $l \times l$ nonnegative symmetric matrix.
\end{theorem}

A remarkable technical achievement of our result is that it does
not require undersmoothing of the nonparametric component, as in
Theorem \ref{thm:ple}, thanks to the profile likelihood approach.
The key lies in a special orthogonality property of the least
favorable curve (see equation (\ref{eqn:defining}) and Lemma
\ref{lemma:order4}). Asymptotic normality without undersmoothing
is also proved in \cite{VanKeilegomCarroll(2007)} for both
backfitting and profiling methods.

Theorem \ref{thm:asymnorm} shows that profile likelihood produces
a semi-parametric efficient estimate even when the number of
parameters diverges. To see this more explicitly, let $p_n = r$ be
a constant. Then, by taking $A_n = I_r$, we obtain
$$
    \sqrt{n}(\hat{\boldsymbol{\beta}}_n -
    \boldsymbol{\beta}_{n0}) \dcv N(0,
    I^{-1}(\boldsymbol{\beta}_{n0})).$$
The asymptotic variance of $\hat{\boldsymbol{\beta}}_n$ achieves
the efficient lower bound given, for example, in
\cite{SeveriniWong(1992)}.

\vspace{12pt}
\subsection{Profile likelihood ratio test}\label{subsect:hyptest}

After estimation of parameters, it is of interest to test the
statistical significance of certain variables in the parametric
component. Consider the problem of testing linear hypotheses:
 $$
    H_0 : A_n\b_0 = 0 \longleftrightarrow  H_1 : A_n\b_{0} \neq 0,
 $$
where $A_n$ is an $l \times p_n$ matrix and $A_nA_n^T = I_l$ for a
fixed $l$. Note that both the null and the alternative hypotheses
are semi-parametric, with nuisance functions $\a(\cdot)$.  The
generalized likelihood ratio test (GLRT) is defined by
 $$  T_n = 2 \{\sup_{\Omega_n}\hat{Q}_n(\b) -
    \sup_{\Omega_n; A_n\b = 0}\hat{Q}_n(\b)\}.
 $$
Note that the testing procedure does not depend explicitly on the
estimated asymptotic covariance matrix. The following theorem
shows that, even when the number of parameters diverges with
sample size, $T_n$ still follows a chi-square distribution
asymptotically, without reference to any nuisance parameters and
functions. This reveals the Wilk's phenomenon, as termed in
\cite{Fanetal(2001)}.

\begin{theorem}\label{thm:hyptest}
Assuming Conditions (A) - (G), under $H_0$, we have
$$T_n \dcv \chi_l^2,$$
provided that $p_n^5/n = o(1)$ and $h=O(n^{-a})$ for $3/(10(p+1))
< a < 2/5$.
\end{theorem}

\subsection{Consistency of the sandwich covariance formula}\label{subsect:sandwichcov}
The estimated covariance matrix for $\bhat$ can be obtained by the
sandwich formula
\begin{equation*}
\hat{\Sigma}_n = n^2\{\nabla^2 \hat{Q}_n(\bhat)\}^{-1}
\widehat{\text{cov}}\{\nabla\hat{Q}_n(\bhat)\}\{\nabla^2
\hat{Q}_n(\bhat)\}^{-1},
\end{equation*}
where the middle matrix has $(j,k)$ entry given by
\begin{equation*}
\begin{split}
(\widehat{\text{cov}}\{\nabla\hat{Q}_n(\bhat)\})_{jk} =
&\biggl\{\frac{1}{n} \sum_{i=1}^n
\frac{\partial\hat{Q}_{ni}(\bhat)}{\partial\beta_{nj}}
\frac{\partial\hat{Q}_{ni}(\bhat)}{\partial\beta_{nk}}\biggr\}\\
   &-\biggl\{\frac{1}{n}\sum_{i=1}^n
\frac{\partial\hat{Q}_{ni}(\bhat)}{\partial\beta_{nj}}
\frac{1}{n}\sum_{i=1}^n
\frac{\partial\hat{Q}_{ni}(\bhat)}{\partial\beta_{nk}}\biggr\}.
\end{split}
\end{equation*}
With the notation $\Sigma_n = I_n^{-1}(\b_0)$, we have the
following consistency result for the sandwich formula.

\begin{theorem}\label{thm:sandwich}
Assuming Conditions (A) - (G). If $p_n^4/n = o(1)$ and
$h=O(n^{-a})$ with $(4(p+1))^{-1} < a < 1/2$, we have
 $$ A_n\hat{\Sigma}_nA_n^T - A_n\Sigma_nA_n^T  \pcv 0
\text{ as } n\rightarrow\infty $$ for any $l \times p_n$  matrix
$A_n$ such that $A_nA_n^T = G$.
\end{theorem}

This result provides a simple way to construct confidence
intervals for $\b$. Simulation results show that this formula
indeed provides a good estimate of the covariance of $\bhat$ for a
variety of practical sample sizes.

\section{Computation of the
estimates}\label{sect:computation}

\setcounter{equation}{0} Finding $\bhat$ to maximize the profile
likelihood (\ref{eqn:glfhat}) poses some interesting challenges,
as the function $\hat{\a}_\b(u)$ in (\ref{eqn:glfhat}) depends on
$\b$ implicitly (except the least-square case). The full
profile-kernel estimate is to directly employ the Newton-Raphson
iterations:
\begin{equation}
\boldsymbol{\beta}_n^{(k+1)} = \boldsymbol{\beta}_n^{(k)} -
\{\nabla^2\hat{Q}_n(\boldsymbol{\beta}_n^{(k)})
    \}^{-1}\nabla \hat{Q}_n(\boldsymbol{\beta}_n^{(k)}),
    \label{NR}
\end{equation}
starting from the initial estimate $\bbeta^{(0)}$.  We will call
the estimate $\boldsymbol{\beta}_n^{(k)}$ and
$\hat{\a}_{\boldsymbol{\beta}_n^{(k)}}(u)$ the $k$-step estimate
(\cite{Bickel(1975)}; \cite{Robinson(1988)}). The initial estimate
for $\b$ is critically important for the computational speed. We
will propose a new and fast initial estimate in Section
\ref{subsect:DBE}.

The first two derivatives of
$\thickspace\thinspace\nabla\hat{Q}_n(\b)$ is given by
\begin{equation}\label{eqn:2.3A}
\begin{split}
&\thickspace\thinspace\nabla\hat{Q}_n(\b) = \sum_{i=1}^n
q_{1i}(\b)
(\mathbf{Z}_{ni} + \hat{\a}_\b^{\prime}(U_i)\mathbf{X}_i),\\
&\nabla^2\hat{Q}_n(\b) = \sum_{i=1}^n q_{2i}(\b)
(\mathbf{Z}_{ni} + \hat{\a}_\b^{\prime}(U_i)\mathbf{X}_i)(\mathbf{Z}_{ni} + \hat{\a}_\b^{\prime}(U_i)\mathbf{X}_i)^T\\
&\;\;\;\;\;\;\;\;\;\;\;\;\;\;\;\;\;\;\;+\sum_{i=1}^n
\biggl\{q_{1i}(\b) \sum_{r=1}^q
\frac{\partial^2\hat{\alpha}_{\b}^{(r)}(U_i)}
{\partial\b\partial\boldsymbol{\beta}_n^T}X_{ir}\biggr\},
\end{split}
\end{equation}
where $q_l(x, y) = \frac{\partial^l}{\partial x^l} Q(g^{-1}(x),
y)$, $q_{ki}(\b)  = q_k(\hat{m}_{ni}(\b), Y_{ni})$ ($k=1, 2$) with
$\hat{m}_{ni}(\b) = \hat{\a}_\b(U_i)^T\mathbf{X}_i +
\mathbf{Z}_{ni}^T\b$. In the above formulae,
$\hat{\a}_\b^{\prime}(u) =
\frac{\partial\hat{\a}_\b(u)}{\partial\b}$ is a $p_n$ by $q$
matrix and $\alpha_{\b}^{(r)}(u)$ is the $r^{\text{th}}$ component
of $\a_{\b}(u)$.

As the first two derivatives of $\hat{\a}_\b(u)$ are hard to
compute in (\ref{eqn:2.3A}), one can employ the backfitting
algorithm, which iterates between (\ref{eqn:llf}) and
(\ref{eqn:glf}). This is really the same as the fully iterated
algorithm (\ref{NR}) but ignores the functional dependence of
$\hat{\a}_\b(u)$ in (\ref{eqn:glfhat}) on $\b$; it uses the value
of $\b$ in the previous step of the iteration as a proxy. More
precisely, the backfitting algorithm treats the terms
$\hat{\a}_\b'(u)$ and $\hat{\a}_\b''(u)$ in (\ref{eqn:2.3A}) as
zero and computes $\hat{m}_{ni}(\b)$ using the value of $\b$ from
the previous iteration. The maximization is thus much easier to
carry out, but the convergence speed can be reduced. See
\cite{Huetal(2004)} and \cite{LinCarroll(2006)} for more
descriptions of the two methods and some closed-form solutions
proposed for the partially linear models.

Between these two extreme choices is our modified algorithm, which
ignores the computation of the second derivative of
$\hat{\a}_\b(u)$ in (\ref{NR}), but keeps its first derivative in
the iteration. Namely, the second term in (\ref{eqn:2.3A}) is
treated as zero. Details will be given in Section
\ref{subsect:betaupdate}. It turns out that this algorithm
improves significantly the computation with achieved accuracy.  At
the same time, it enhances dramatically the stability of the
algorithm.  We will term the algorithm as the accelerated
profile-kernel algorithm.

When the quasi-likelihood becomes a square loss, the accelerated
profile-kernel algorithm is exactly the same as that used to
compute the full profile likelihood estimate, since $\ahat(\cdot)$
is linear in $\b$.

\vspace{12pt}
\subsection{Difference-based estimation}\label{subsect:DBE}

We generalize the difference-based idea to obtain an initial
estimate $\boldsymbol{\beta}_n^{(0)}$. The idea has been used in
\cite{Yatchew(1997)} and \cite{FanHuang(2005)} to remove the
nonparametric component in the partially linear model.

We first consider the specific case of the GVCPLM:
\begin{equation}\label{eqn:VCPLM}
Y = \a(U)^T\mathbf{X} + \b^T\mathbf{Z}_n + \varepsilon.
\end{equation}
This is the varying-coefficient partially linear model studied by
\cite{Zhangetal(2002)} and \cite{Xiaetal(2004)}. Let the random
sample $\{(U_i,\mathbf{X}_i^T, \mathbf{Z}_{ni}^T, Y_i)\}_{i=1}^n$
be from the model (\ref{eqn:VCPLM}), with the data ordered
according to the $U_i$'s. Under mild conditions, the spacing
$U_{i+j} - U_i$ is $O_P(1/n),$ so that
\begin{equation}
    \a(U_{i+j}) - \a(U_i) \approx \boldsymbol{\gamma_0} +
\boldsymbol{\gamma_1}(U_{i+j}-U_i), \qquad j=1,\cdots,q.
\label{approx}
\end{equation}
Indeed, it can be approximately zero; the linear term is used to
reduce the approximation errors.

For given weights $w_j$ (its dependence on $i$ is suppressed for
simplicity), define
$$
Y_i^* = \sum_{j=1}^{q+1}w_jY_{i+j-1}, \quad \mathbf{Z}_{ni}^* =
\sum_{j=1}^{q+1}w_j\mathbf{Z}_{n(i+j-1)}, \quad \varepsilon_i^* =
\sum_{j=1}^{q+1}w_j\varepsilon_{i+j-1}.
$$
If we choose the weights to satisfy
$\sum_{j=1}^{q+1}w_j\mathbf{X}_{i+j-1} = \mathbf{0}$, then using
(\ref{eqn:VCPLM}) and (\ref{approx}), we have
\begin{equation*}
Y_i^* \approx \boldsymbol{\gamma_0}^T\mathbf{X}_iw_1 +
\boldsymbol{\gamma_1}^T\sum_{j=1}^{q+1}w_j
U_{i+j-1}\mathbf{X}_{i+j-1}+
\boldsymbol{\beta}_n^T\mathbf{Z}_{ni}^* + \varepsilon_i^*,
\end{equation*}
Ignoring the approximation, which is of order $O_P(n^{-1})$, the
above is a multiple regression model with parameters
$(\boldsymbol{\gamma_0}, \boldsymbol{\gamma_1}, \b)$. The
parameters can be found by a weighted least square fit to the
$(n-q)$ starred data. This yields a root-n consistent estimate of
$\b$, as the above approximation for the finite $q$ is of order
$O_P(n^{-1})$.

To solve $\sum_{j=1}^{q+1}w_j\mathbf{X}_{i+j-1}=\mathbf{0}$, we
need to find the rank of the matrix
$(\mathbf{X}_i,\cdots,\mathbf{X}_{i+q})$, denoted it by $r$. Fix
$q+1-r$ of the $w_j$'s and the rest can be determined uniquely by
solving the system of linear equations for $\{w_j, j=1, \cdots,
q+1\}$.  For random designs, with probability 1, $r=q$.  Hence,
the direction of the weights $\{w_j, j=1, \cdots, q+1\}$ is
uniquely determined.  For example, in the partial linear model,
$q=1$ and $\bX_i = 1$.  Hence, $(w_1, w_2) = c(1, -1)$ and the
constant $c$ can be taken to have a norm one.  This results in the
difference based estimator in \cite{Yatchew(1997)} and
\cite{FanHuang(2005)}.

To use the differencing idea to obtain an initial estimate of $\b$
for the GVCPLM, we apply the transformation of the data. If $g$ is
the link function, we use $g(Y_i)$ as the transformed data and
proceed with the difference-based method as for the VCPLM. Note
that for some models like the logistic regression with logit link
and Poisson log-linear model, we need to make adjustments in
transforming the data. We use $g(y) =
\log(\frac{y+\delta}{1-y+\delta})$ for the logistic regression and
$g(y) = \log(y+\delta)$ for the Poisson regression. Here, the
parameter $\delta$ is treated as a smoothing parameter like $h$,
and its choice will be discussed in Section
\ref{subsect:bandwidth}.

\subsection{Accelerated profile-kernel algorithm} \label{subsect:betaupdate}

As mentioned before, the accelerated profile-kernel algorithm
needs to compute $\a_{\b}^{\prime}(u)$, which will be replaced by
its consistent estimate given in the following theorem. The proof
is in section \ref{sect:proof}.

\begin{theorem}\label{thm:consistentalphaprime}
Under Regularity Conditions (A)-(G), provided $\sqrt{p_n}(h +
c_n\lh)=o(1)$ where $c_n = (nh)^{-1/2}$, we have for each
$\b\in\Omega_n$,
$$
    \hat{\a}_\b^{\prime}(u) = - \biggl\{\sum_{i=1}^n
    q_{2i}(\b)
    \mathbf{Z}_{ni}\mathbf{X}_i^TK_h(U_i-u)\biggr\}
\cdot\biggl\{\sum_{i=1}^n q_{2i}
(\b)\mathbf{X}_i\mathbf{X}_i^TK_h(U_i-u) \biggr\}^{-1}
$$
being a consistent estimator of $\a_{\b}^{\prime}(u)$ which holds
uniformly in $u\in\Omega$.
\end{theorem}

Since the function $q_2(\cdot,\cdot)<0$ by Regularity Condition
(D), by ignoring the second term in (\ref{eqn:2.3A}), the modified
$\nabla^2\hat{Q}_n(\b)$ in equation (\ref{eqn:2.3A}) is still
negative-definite. This ensures the Newton-Raphson update of the
profile-kernel procedure can be carried out smoothly. The
intuition behind the modification is that, for a neighborhood
around the true parameter $\b_0$, the least favorable curve
$\a_{\b}(u)$ should be approximately linear in $\b$.

\subsection{One-step estimation for the nonparametric component}
\label{subsect:onestep}

Given $\b = \boldsymbol{\beta}_n^{(k)}$, we need to compute
$\a_\b(u)$ in order to compute $\hat{m}_{ni}(\b)$ and hence the
modified gradient vector and Hessian matrix in (\ref{NR}).  This
is the same as estimating the varying coefficient functions under
model (\ref{eqn:GVCPLM}) with known $\b$.  \cite{Caietal(2000)}
propose a one-step local MLE, which is shown to be as efficient as
the fully iterated one.  They also propose an efficient algorithm
to compute these varying coefficient functions.  Their algorithm
can be directly adapted here.  Details can be found in
\cite{Caietal(2000)}.

\subsection{Choice of bandwidth}\label{subsect:bandwidth}

As mentioned at the end of Section \ref{subsect:DBE}, in addition
to choosing the bandwidth $h$, we have an extra smoothing
parameter $\delta$ to be determined due to the adjustments to the
transformation of the response $Y_{ni}$. This two dimensional
smoothing parameters $(\delta, h)$ can be selected by a $K$-fold
cross-validation, using the quasi-likelihood as a criterion
function. As demonstrated in Section \ref{sect:simulation}, the
practical accuracy can be achieved in several iterations using the
accelerated profile-kernel algorithm.  Hence, the profile-kernel
estimate can be computed rapidly. As a result, the $K$-fold
cross-validation is not too computationally intensive, as long as
$K$ is not too large (e.g. K=5 or 10).

\section{Numerical properties}\label{sect:simulation}

\setcounter{equation}{0} To evaluate the performance of estimator
$\hat{\a}(\cdot)$, we use the square-root of average errors (RASE)
$$\text{RASE} = \biggl\{n_{\text{grid}}^{-1}\sum_{k=1}^{n_{\text{grid}}}\|\hat{\a}(u_k)-\a(u_k)\|^2 \biggr\}^{1/2}, $$
over $n_{\text{grid}} = 200$ grid points $\{u_k\}$.  The
performance of the estimator $\bhat$ is assessed by the
generalized mean square error (GMSE)
$$ \text{GMSE} = (\bhat - \b_0)^T \bB (\bhat-\b_0),$$
where $\bB = E \mathbf{Z}_n \mathbf{Z}_n^T $.

Throughout our simulation studies, the dimensionality of
parametric component is taken as $p_n = \lfloor 1.8n^{1/3}\rfloor$
and the nonparametric component as $q = 2$ in which $X_1 = 1$ and
$X_2 \sim N(0, 1)$. The rate $p_n = O_P(n^{1/3})$ is not the same
as presented in the theorems in section \ref{sect:properties}, but
we use this to show the capability of handling a higher rate of
parameters growth for the accelerated profile-kernel method. In
addition, the covariates $(\mathbf{Z}_n^T, X_2)^T$ is a
$(p_n+1)-$dimensional normal random vector with mean zero and
covariance matrix $(\sigma_{ij})$, where $\sigma_{ij} =
0.5^{|i-j|}$.  Furthermore, we always take $U \sim U(0,1)$
independent of the other covariates.  Finally, we use
$\text{SD}_{\text{mad}}$ to denote the robust estimate of standard
deviation, which is defined as interquartile range divided by
1.349. The number of simulations is 400 except that in Table 1
(which is 50) due to the intensive computation of the fully
iterated profile-kernel estimate.

\noindent\textbf{Poisson model}. The response $Y$, given
$(U,\mathbf{X},\mathbf{Z_n})$, has a Poisson distribution with the
mean function $\mu(U,\mathbf{X},\mathbf{Z}_n)$ where
  $$
    \log(\mu(U,\mathbf{X},\mathbf{Z}_n)) =
      \mathbf{X}^T\a(U) + \mathbf{Z}_n^T\b.
  $$
We have $\b_0 = (0.5, 0.3, -0.5, 1, 0.1, -0.25,0,\cdots,0)^T$, the
$p_n$-dimensional parameters. The coefficient functions are given
by
$$
   \alpha_1(u) = 4 + \sin(2\pi u), \;\;\text{and }\; \alpha_2(u) = 2u(1-u).
$$

\noindent\textbf{Bernoulli model.}  The response $Y$, given
$(U,\mathbf{X},\mathbf{Z_n})$, has a Bernoulli distribution with
the success probability given by
$$
    p(U,\mathbf{X},\mathbf{Z}_n)) = \exp\{\mathbf{X}^T\a(U) + \mathbf{Z}_n^T\b\}
    /[1+\exp\{\mathbf{X}^T\a(U) + \mathbf{Z}_n^T\b\}].
$$
The $p_n-$dimensional parameters are $\b_0 = (3, 1, -2, 0.5, 2,
-2,0,\cdots,0)^T$ and the varying coefficient functions is given
by
$$\alpha_1(u) = 2(u^3+2u^2-2u), \;\;\text{and }\;
\alpha_2(u) = 2\cos(2\pi u). $$

Throughout our numerical studies, we use the Epanechnikov kernel
$K(u) = 0.75(1-u^2)_+$ and the 5-fold cross-validation to choose a
bandwidth $h$ and $\delta$. With the assistance of the 5-fold
cross-validation, we chose $\delta = 0.1$ and $h=0.1, 0.08, 0.075$
and $0.06$ respectively for $n=200, 400, 800$ and 1500 for the
Poisson model. For the Bernoulli model, $\delta = 0.005$ and $h =
0.45, 0.4, 0.25$ and 0.18 were chosen respectively for
$n=200,400,800$ and 1500.

Note that $X_2$ and the $Z_{ni}$'s are not bounded r.v.s as needed
in condition (A) in section \ref{sect:proof}. However, these still
satisfy the moment conditions needed in the proofs, and condition
(A) is imposed to merely simplify these proofs.  Condition (B) is
satisfied mainly because the correlations between further
$Z_{ni}$'s are weak, and condition (C) is satisfied because it
involves products of standard normal r.v.s which are bounded in
the first  two moments.

\subsection{Comparisons of algorithms}

\begin{table}[htbp]
  \centering
  \caption{Computation time and accuracy for different
   computing algorithms}
  \begin{tabular}{cc|ccc}\label{table:1}\\
  \hline
    $n$ & $p_n$     & backfitting  & accelerated profile-kernel & full profile-kernel\\
  \hline
  &    \multicolumn{4}{c}{\emph{Median and $\text{SD}_\text{mad}$ (in parentheses)
        of computing times in seconds}}\\
  \hline
  200 & 10     & .6(.0)    &  .7(.0)    &  77.2(.2)  \\
  400 & 13     & .8(.0)    &  1.4(.0)    & 463.2(.9)  \\
  \hline
    &   \multicolumn{4}{c}{\emph{Median and $\text{SD}_\text{mad}$
    (in parentheses) of GMSE (multiplied by $10^4$)}}\\
  \hline

  200 & 10     & 10.72(6.47)   & 5.45(2.71)      & 9.74(14.67)    \\
  400 & 13     & 5.63(4.39)    & 2.78(1.19)      & 5.26(9.46)     \\
  \hline
    &  \multicolumn{4}{c}{\emph{Median RASE relative to the oracle estimate}}\\
  \hline
  200 & 10     & .848    &  .970    &  .895  \\
  400 & 13     & .856    &  .986    &  .882  \\
  \hline
\end{tabular}
\end{table}

We first compare the computing times and the accuracies among
three algorithms: 3-step backfitting, 3-step accelerated
profile-kernel and fully-iterated profile-kernel algorithms. All
of them use the difference-based estimate as the initial estimate.
Table \ref{table:1} summarizes the results based on the Poisson
model with 50 samples.

With the same initial values, the backfitting algorithm is
slightly faster than the accelerated profile-kernel algorithm,
which in turn by far faster than the full profile-kernel
algorithm.  Our experience shows that the backfitting algorithm
needs more than 20 iterations to converge without improving too
much the GMSE.  In terms of the accuracy of estimating the
parametric component, the accelerated profile-kernel algorithm is
about twice as accurate as the backfitting algorithm and the full
profile-kernel one.  This demonstrates the advantage of keeping
the curvature of the least-favorable function in the
Newton-Raphson algorithm.  For the nonparametric component, we
compare RASEs of three algorithms with those based on the oracle
estimator, which uses the true value of $\b$.  The ratios of the
RASEs based on the oracle estimator and those based on the three
algorithms are reported in Table 1.  It is clear that the
accelerated profile-kernel estimate performs very well in
estimating the nonparametric components, mimicking very well the
oracle estimator. The second best is the backfitting algorithm.

We have also compared the three algorithms using the Bernoulli
model. Our proposed accelerated profile-kernel estimate still
performs the best in terms of accuracy, though the improvement is
not as dramatic as those for the Poisson model.  We speculate that
the poor performance of the full profile-kernel estimate is due to
its unstable implementation that is related to computing the
second derivatives of the least-favorable curve.

\begin{table}[htbp]
  \centering
  \caption{Medians of the percentages of GMSE
  based on the accelerated profile-kernel estimates}
  \begin{tabular}{cc|cc|cc}\label{table:2}\\
  \hline
       &       & \multicolumn{2}{c|}{\emph{Poisson}} & \multicolumn{2}{c}{\emph{Bernoulli}} \\
   $n$ & $p_n$      & AF/DBE & AF/3S      & AF/DBE & AF/3S        \\
  \hline
  200 & 10          & 8.2    & 99.9       & 64.1   & 101.7        \\
  400 & 13          & 6.0    & 100.2      & 52.7   & 104.7      \\
  800 & 16          & 5.0    & 100.1      & 50.9   & 102.6      \\
  1500 & 20         & 4.2    & 100.0      & 46.4   & 100.5      \\
  \hline
\end{tabular}
\end{table}

We next demonstrate the accuracy of the three-step accelerated
profile-kernel estimate (3S), compared with the fully-iterated
accelerated profile-kernel estimate (AF) (iterating until
convergence), and the difference-based estimate (DBE), which is
our initial estimate. Table \ref{table:2} reports the ratios of
GMSE based on 400 simulations. It demonstrates convincingly that
with the DBE as the initial estimate, three iterations achieve the
accuracy that is comparably with the fully iterated algorithm. In
fact, the one-step accelerated profile-kernel estimates improve
dramatically (not shown here) our initial estimate (DBE). On the
other hand, the DBE itself is not accurate enough for GCVPLM.

\begin{table}
  \centering
  \caption{One-step estimate of parametric components with different bandwidths}
  \begin{tabular}{cc|cc|cc|cc}\label{table:3}\\
  \hline
       &       & \multicolumn{4}{c|}{\emph{Poisson}} & \multicolumn{2}{c}{\emph{Bernoulli}}\\
  \hline
       &       & \multicolumn{2}{c|}{\emph{Median and $\text{SD}_\text{mad}$ of }}&
       \multicolumn{2}{c|}{\emph{Mean and SD of }} &
       \multicolumn{2}{c}{\emph{Median and $\text{SD}_\text{mad}$ of }}\\
       & &
       \multicolumn{2}{c|}{\emph{GMSE}$\times 10^5$}
       &\multicolumn{2}{c|}{\emph{MSE $\times 10^4$ for $\beta_5$}}
       &
       \multicolumn{2}{c}{\emph{GMSE}$\times 10$}\\
   $n$ & $p_n$ & $h_{\text{CV}}$ & $1.5h_{\text{CV}}$ & $0.66h_{\text{CV}}$ & $h_{\text{CV}}$ & $0.66h_{\text{CV}}$ & $h_{\text{CV}}$\\
  \hline
   200 &  10   &5.9(3.0)     &6.4(3.3)     &993(112)  &995(105)   &8.2(4.4)     &8.4(5.1)    \\
   400 &  13   &3.1(1.4)     &3.0(1.4)     &1004(67)  &1001(65)   &4.8(2.2)     &5.4(2.5)    \\
   800 &  16   &1.7(0.7)     &1.7(0.6)     &999(47)   &999(46)    &2.7(1.0)     &2.7(1.1)    \\
   1500&  20   &1.1(0.3)     &1.1(0.4)     &1000(32)  &1000(32)   &1.8(0.7)     &1.8(0.6)    \\
   \hline
  \end{tabular}\\
  {\scriptsize SD and $\text{SD}_\text{mad}$ are shown in
  parentheses.}
\end{table}

The effect of bandwidth choice on the estimation of parametric
component is summarized  in Table \ref{table:3}. Denote by
$h_\text{CV}$ the bandwidth chosen by the cross-validation.  We
scaled the bandwidth up and down by using a factor of 1.5.  For
illustration, we use the one-step accelerated profile-kernel
estimate. The results for three-step profile-kernel estimate are
similar. We evaluate the performance for all components using GMSE
and for the specific component $\beta_5$ using MSE (the results
for other components are similar). We do not report all the
results here to save the space. It is clear that the GMSE does not
sensitively depends on the bandwidth, as long as it is reasonably
close to $h_\text{CV}$.  This is consistent with our asymptotic
results.

\subsection{Accuracy of profile-likelihood inferences}

\begin{table}[htbp]
  \centering
  \caption{Standard deviations and estimated standard errors}
  \begin{tabular}{cc|cc|cc|cc|cc}\label{table:4}\\
  \hline
       &    & \multicolumn{4}{c|}{\emph{Poisson, values$\times$1000}} & \multicolumn{4}{c}{\emph{Bernoulli, values$\times$10}}\\
       &    & \multicolumn{2}{c}{$\hat{\beta}_1$}& \multicolumn{2}{c|}{$\hat{\beta}_3$}& \multicolumn{2}{c}{$\hat{\beta}_2$} &\multicolumn{2}{c}{$\hat{\beta}_4$}\\
\hline
   $n$ & $p_n$ & SD & $\text{SD}_{\text{m}}$ & SD &$\text{SD}_{\text{m}}$& SD &$\text{SD}_{\text{m}}$& SD &
   $\text{SD}_{\text{m}}$ \\
  \hline
   200 & 10 &9.1&8.5(1.3)&9.9&9.4(1.3)&3.6&2.9(.4)&3.2&2.8(.4)\\
   400 & 13 &6.0&5.6(0.7)&6.5&6.1(0.7)&2.3&2.1(.2)&2.2&2.0(.2)\\
   800 & 16 &3.7&3.8(0.3)&4.1&4.2(0.4)&1.7&1.6(.1)&1.5&1.5(.1)\\
   1500& 20 &2.8&2.7(0.2)&3.1&3.0(0.2)&1.2&1.2(.1)&1.1&1.1(.1)\\
   \hline
  \end{tabular}\\
  {\scriptsize $\text{SD}_\text{mad}$ are shown in
  parentheses.}
\end{table}

To test the accuracy of the sandwich formula for estimating
standard errors, the standard deviations of the estimated
coefficients (using the one-step accelerated profile-kernel
estimate) are computed from the 400 simulations using
$h_{\text{CV}}$. These can be regarded as the true standard errors
(columns labeled SD). The 400 estimated standard errors are
summarized by their median (columns $\text{SD}_m$) and its
associated $\text{SD}_{\text{mad}}$. Table 4 summarizes the
results. Clearly, the sandwich formula does a good job, and
accuracy gets better as $n$ increases.


We now study the performance of GLRT in Section
\ref{subsect:hyptest}. To this end, we consider the following null
hypothesis:
$$ H_0:\beta_7=\beta_8=\cdots=\beta_{p_n} = 0.$$
We examine the power of the test under a sequence of the
alternative hypotheses indexed by a parameter $\gamma$ as follows:
$$ H_1:\beta_7=\beta_8=\gamma, \; \beta_j=0 \text{ for } j>8.$$
When $\gamma=0$, the alternative hypothesis becomes the null
hypothesis.

\begin{figure}[htbp]
\centerline{\psfig{figure=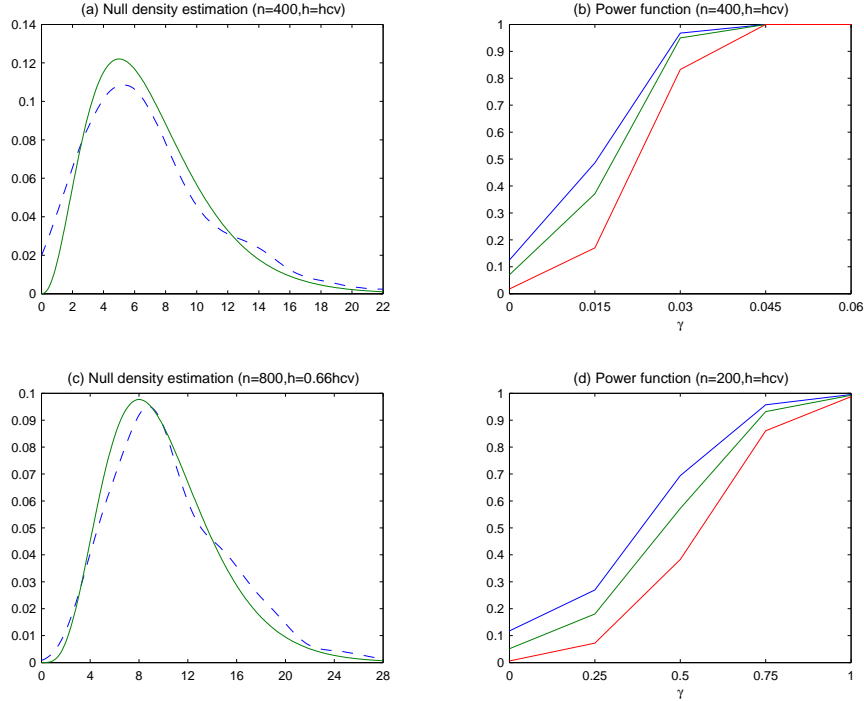,width=4.5in}}
\caption[Fig1] {\emph{(a) Asymptotic null distribution (solid) and
estimated true null distribution (dotted) for the Poisson model.
(b) The power function at significant level $\alpha = 0.01, 0.05$
and 0.1.  The captions for (c) and (d) are the same as those in
(a) and (b) except that the Bernoulli model is now used.}}
\label{figure:1}
\end{figure}

Under the null hypothesis, the GLRT statistics are computed for
each of 400 simulations, using the one-step accelerated
profile-kernel estimates. Their distribution is summarized by a
kernel density estimate and can be regarded as the true null
distribution. This is compared with the asymptotic null
distribution $\chi^2_{p_n - 6}$. Figures \ref{figure:1}(a) and (c)
show the results when $n=400$. The finite sample null density is
seen to be reasonably close to the asymptotic one, except for the
Monte Carlo error.

The power of the GLR test is studied under a sequence of
alternative models, progressively deviating from the null
hypothesis, namely, as $\gamma$ increases.  Again, the one-step
accelerated profile-kernel algorithm is employed. The power
functions are calculated at three significance levels: 0.1, 0.05
and 0.01, using the asymptotic distribution. They are the
proportion of rejection among the 400 simulations and are depicted
in Figures \ref{figure:1}(b) and (d). The power curves increase
rapidly with $\gamma$, which shows the GLR test is powerful.  The
powers at $\gamma = 0$ are approximately the same as the
significance level except the Monte Carlo error.  This shows that
the size of the test is reasonably accurate.

\subsection{A real data example.}\label{subsect:realdataanalysis}

This is the analysis of the data in section \ref{subsect:example}
in where details of data and variables are given.

To examine the nonlinear effect of age and its nonlinear
interaction with the experience, we appeal to the following GVCPLM
(interactions between age and covariates other than
\textbf{TotalYrsExp} are considered but found to be
insignificant):

\begin{equation}\label{eqn:4.1}
\begin{split}
\log\bigg(\frac{p_H}{1-p_H}\bigg) =
& \:\alpha_1(\text{Age}) + \alpha_2(\text{Age})\text{TotalYrsExp}\\
&+ \beta_1\text{Female} + \beta_2\text{PCJob} +
\sum_{i=1}^4\beta_{2+i}\text{Edu}_i
\end{split}
\end{equation}
where $p_H$ is the probability of having a high grade job.
Formally, we are testing
\begin{equation} \label{eqn:4.2}
 H_0: \beta_1 = 0 \longleftrightarrow
H_1: \beta_1 < 0.
\end{equation}

\begin{table}[htbp]
  \centering
  \caption{Fitted coefficients (sandwich SD) for model (\ref{eqn:4.1})}
  \begin{tabular}{ccccccc}\label{table:5}\\
  \hline
    \emph{Response}  & Female & PCJob & $\text{Edu}_1$ & $\text{Edu}_2$ & $\text{Edu}_3$ & $\text{Edu}_4$   \\
  \hline
  HighGrade4 & -1.96(.57) &  -0.02(.76) &  -5.14(.85) &  -4.77(.98) &  -2.72(.52) &  -2.85(.96) \\
  HighGrade5 & -2.22(.59) &  -1.96(.61) &  -5.69(.67) &  -5.95(.97) &  -3.09(.72) &  -1.26(1.10) \\
  \hline
\end{tabular}
\end{table}

\begin{figure}[htbp]
\centerline{\psfig{figure=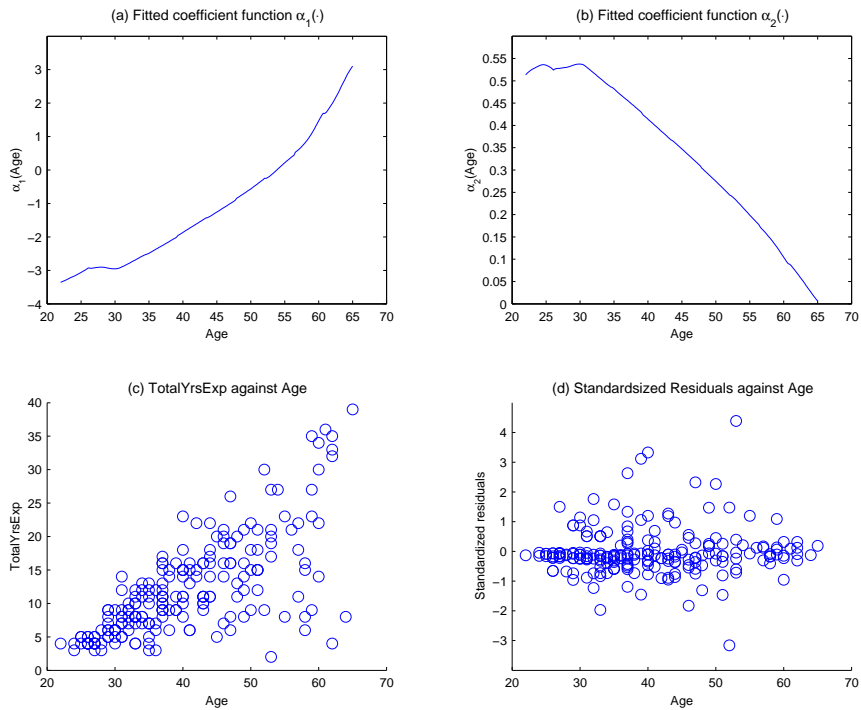,width=4.5in}}
\caption
{\emph{(a) Fitted coefficient function $\alpha_1(\cdot)$ (b)
 Fitted coefficient function $\alpha_2(\cdot)$.
 (c) The scatter plot `TotalYrsExp' Against `Age'. (d) Standardized residuals against the
variable `Age'. }}\label{Fig2}
\end{figure}

A 20-fold CV is employed to select the bandwidth $h$ and the
parameter $\delta$ in the transformation of the data.  This yields
$h_\text{CV}=24.2$, $\delta_\text{CV} = 0.1$.   Table
\ref{table:5} shows the results of the fit using the three-step
accelerated profile-kernel estimate. The coefficient for
\textbf{Female} is significantly negative. The education plays
also an important role in getting high grade job. All coefficients
are negative, as they are contrasted with the highest education
level. The \textbf{PCJob} does not seem to play any significant
role in getting promotion. Figures \ref{Fig2}(a) and (b) depict
the estimated coefficient functions. They show that as age
increases one has a better chance of being in a higher job grade,
and then the marginal effect of working experience is large when
age is around 30 or less, but start to fall as one gets older.
However, the second result should be interpreted with caution, as
the variables \textbf{Age} and \textbf{TotalYrsExp} are highly
correlated (Figure \ref{Fig2}(c)).  The standardized residuals
$(y-\hat{p}_\text{H})/\sqrt{\hat{p}_\text{H}(1-\hat{p}_\text{H})}$
against \textbf{Age} is plotted in Figure \ref{Fig2}(d).  It shows
that the fit seems reasonable. Other diagnostic plots also look
reasonable, but they are not shown here.

We have conducted another fit using a binary variable
\textbf{HighGrade5}, which is 0 only when job grade is less than
5. The coefficients are shown in Table \ref{table:5} and the
\textbf{Female} coefficient is close to the first fit.

We now employ the generalized likelihood ratio test to the problem
(\ref{eqn:4.2}).  The GLR test statistic is 14.47 with one degree
of freedom, resulting in a P-value of 0.0001. We have also conduct
the same analysis using \textbf{HighGrade5} as the binary
response.  The GLR test statistic is now 13.76 and the associated
P-value is 0.0002.  The fitted coefficients are summarized in
Table 5.  The result provides stark evidence that even after
adjusting for other confounding factors and variables, female
employees  of the Fifth National Bank is harder to get promoted to
a high grade job.

Not shown in this paper, we have conducted the analysis again
after deleting 6 data points corresponding to 5 male executives
and 1 female employee having many years of working experience and
high salaries. The test results are still similar.

\section{Technical proofs}\label{sect:proof}

\setcounter{equation}{0} In this section the proofs of Theorems
1-4 will be given. We introduce some notations and regularity
conditions for our results. In the following and thereafter, the
symbol $\otimes$ represents the Kronecker product between
matrices, and $\lambda_{\text{min}}(A)$ and
$\lambda_{\text{max}}(A)$ denotes respectively the minimum and
maximum eigenvalues of a symmetric matrix A. We let $Q_{ni}(\b)$
be the $i$-th summand of (\ref{eqn:glf}).

Denote the true linear parameter by $\boldsymbol{\beta}_{n0}$,
with parameter space $\Omega_n \subset \mathbb{R}^{p_n}$. Let
$\mu_k = \int_{-\infty}^{\infty} u^k K(u) du$ and $A_p(\mathbf{X})
= (\mu_{i+j})_{0\leq i,j \leq p} \otimes \X\X^T$.  Set
\begin{eqnarray*}
 \rho_l(t) & = (dg^{-1}(t)/dt)^l/V(g^{-1}(t)), \;\;\; &
m_{ni}(\boldsymbol{\beta}_n)=
\boldsymbol{\alpha}_{\boldsymbol{\beta}_n}(U_i)^{T}\mathbf{X}_i +
\boldsymbol{\beta}_n^{T}\mathbf{Z}_{ni},\\
 & \a_{\b}^{\prime}(u) = \frac{\partial\a_{\b}(u)}{\partial\b} , &
\a_{\b}^{(r)\prime\prime}(u) =
\frac{\partial^2\a_\b^{(r)}(u)}{\partial\b\partial\boldsymbol{\beta}_n^T}.
\end{eqnarray*}

\textbf{Regularity Conditions:}
\begin{itemize}
\item[(A)] The covariates $\mathbf{Z}_{n}$ and $\X$ are bounded
random variables.

\item[(B)] The smallest and the largest eigenvalues of the matrix
$I_n(\boldsymbol{\beta}_{n0})$ is bounded away from zero and
infinity for all $n$.  In addition,
$\mathbf{E}_0[\nabla^T{Q}_{n1}(\boldsymbol{\beta}_{n0})\nabla{Q}_{n1}(\boldsymbol{\beta}_{n0})
]^4 = O(p_n^4)$.

\item[(C)]
$\mathbf{E}_{\boldsymbol{\beta}_n}|\frac{\partial^{l+j}Q_{n1}(\b)}
{\partial^j\a\partial\beta_{nk_1}\cdots\partial\beta_{nk_l}}|$ and
$\mathbf{E}_{\boldsymbol{\beta}_n}|\frac{\partial^{l+j}Q_{n1}(\b)}
{\partial^j\a\partial\beta_{nk_1}\cdots\partial\beta_{nk_l}}|^2$
are bounded for all $n$, with $l = 1,\cdots,4$ and $j=0,1$.

\item[(D)] The function $q_2(x,y) < 0$ for $x\in\mathbb{R}$ and
$y$ in the range of the response variable, and $
\mathbf{E}_0\{q_2(m_{n1}(\boldsymbol{\beta}_n),Y_{n1})A_p(\mathbf{X_1})|U=u
\}$ is invertible.

\item[(E)]
 The functions $V^{\prime\prime}(\cdot)$ and
$g^{\prime\prime\prime}(\cdot)$  are continuous. The
least-favorable curve
$\boldsymbol{\alpha}_{\boldsymbol{\beta}_n}(u)$ is three times
continuously differentiable in $\boldsymbol{\beta}_n$ and $u$.

\item[(F)] The random variable $U$ has a compact support $\Omega$.
The density function $f_U(u)$ of $U$ has a continuous second
derivative and is uniformly bounded away from zero.

\item[(G)] The kernel K is a bounded symmetric density function
with bounded support.

\end{itemize}

Note the above conditions are assumed to hold uniformly in
$u\in\Omega$. Condition (A) is imposed just for the simplicity of
proofs. The boundedness of covariates is imposed to ensure various
products involving $q_l(\cdot,\cdot), \mathbf{X}$ and
$\mathbf{Z}_n$ have bounded first and second moments. Conditions
(B) and (C) are uniformity conditions on higher-order moments of
the likelihood functions. They are stronger than those of the
usual asymptotic likelihood theory, but they facilitate technical
proofs. Condition (G) is also imposed for simplicity of technical
arguments. All of these conditions can be relaxed at the expense
of longer proofs.

Before proving Theorem \ref{thm:ple}, we need two important
lemmas. Lemma \ref{lemma:order1} concerns the order approximations
to the least-favorable curve $\a_\b(\cdot)$, while Lemma
\ref{lemma:order4} holds the key to showing why undersmoothing is
not needed in Theorems \ref{thm:ple} and \ref{thm:asymnorm}. Let
$c_n = (nh)^{-1/2}$, $\hat{\mathbf{a}}_{0\b}$, $\cdots$, and
$\hat{\mathbf{a}}_{p\b}$ maximize (\ref{eqn:llf}), and
$\a_{u\b}^{(p)}(u) = \frac{\partial^p\a_\b(u)}{\partial{u^p}}$.
Set
\begin{eqnarray*}
  \abar(u) &=& \X_i^T\biggl( \sum_{k=0}^p \frac{(U_i-u)^k}{k!}\a_{u\b}^{(k)}(u)   \biggr) + \boldsymbol{\beta}_n^T\Z_{ni}, \\
  \bhatstar &=&  c_n^{-1}\biggl((\hat{\mathbf{a}}_{0\b} - \a_\b(u))^T, \cdots, \frac{h^p}{p!}(\hat{\mathbf{a}}_{p\b}-\a_{u\b}^{(p)}(u))^T \biggr)^T,\\
  \X_i^*     &=& \biggl(1, \frac{U_i-u}{h}, \cdots,
  \biggl(\frac{U_i-u}{h}\biggr)^p\biggr)^T \otimes\X_i.
\end{eqnarray*}

\begin{lemma}\label{lemma:order1} Under Regularity Conditions (A)
- (G), for each $\b\in\Omega_n$, the following holds uniformly in
$u \in \Omega$:
$$ \|\hat{\mathbf{a}}_{0\b}(u) - \a_\b(u) \| = O_P(h^{p+1} +
c_n\lh).$$ Likewise, the norm of the $k^\text{th}$ derivative of
the above with respect to any $\beta_{nj}$'s, for $k=1,\cdots,4$,
all have the same order uniformly in $u \in \Omega$.
\end{lemma}

\noindent {\bf Proof of Lemma \ref{lemma:order1}}. Our first step
is to show that, uniform in $u\in\Omega$,
 $$\bhatstar = \tilde{\mathbf{A}}_n^{-1}\mathbf{W}_n +
O_P(h^{p+1} + c_n\log^{1/2}(1/h)),$$where
\begin{eqnarray*}
  \tilde{\mathbf{A}}_n &=& f_U(u)E_0\{ \rho_2(\a_\b(U)^T\X + \Z_n^T\b)A_p(\X)|U=u\}, \\
  \mathbf{W}_n  &=& hc_n\sum_{i=1}^n q_1(\abar,
  Y_{ni})\X_i^*K_h(U_i-u),\\
  \mathbf{A}_n &=& hc_n^2\sum_{i=1}^n q_2(\abar,
  Y_{ni})\X_i^*\X_i^{*T}K_h(U_i-u).
\end{eqnarray*}
Since expression (\ref{eqn:llf}) is maximized at
$(\hat{\mathbf{a}}_{0\b}, \cdots,  \hat{\mathbf{a}}_{p\b})^T$,
$\bhatstar$ maximizes
\begin{equation*}
\begin{split}
l_n(\bstar) &= h\sum_{i=1}^n\{ Q(g^{-1}(c_n\X_i^{*T}\bstar +
\abar), Y_{ni}) - Q(g^{-1}(\abar), Y_{ni}) \}\\
&= \mathbf{W}_n^T\bstar +
\frac{1}{2}\boldsymbol{\beta}^{*T}\mathbf{A}_n\bstar +
\frac{hc_n^3}{6}\sum_{i=1}^n q_3(\eta_i,
y_{ni})(\X_i^{*T}\bstar)^3K_h(U_i-u),
\end{split}
\end{equation*}
where $\eta_i$ lies between $\abar$ and $\abar +
c_n\X_i^{*T}\bstar $. The concavity of $l_n(\bstar)$ is ensured by
Condition (D). Note that $K(\cdot)$ is bounded, we have under
Conditions (A) and (C), the third term on the right hand side is
bounded by
$$ O_P(nhc_n^3E|q_3(\eta_1, Y_{n1})\|\X_1\|^3K_h(U_1-u)|) =
O_P(c_n) = o_P(1).$$ Direct calculation yields $E_0\mathbf{A}_n =
-\tilde{\mathbf{A}}_n + O(h^{p+1})$ and
$\text{Var}_0((\mathbf{A}_n)_{ij}) = O((nh)^{-1})$ so that
mean-variance decomposition yields
$$
\mathbf{A}_n = -\tilde{\mathbf{A}}_n + O_P(h^{p+1}).
$$
Hence we have
\begin{equation}\label{eqn:maxeqn}
    l_n(\bstar) = \mathbf{W}_n^T\bstar -
    \frac{1}{2}\boldsymbol{\beta}^{*T}\tilde{\mathbf{A}}_n\bstar + o_P(1).
\end{equation}

Note that $\mathbf{A}_n$ is a sum of i.i.d. random variables of
kernel form, by a result of \cite{MackSilverman(1982)},
\begin{equation}\label{eqn:An}
\mathbf{A}_n = -\tilde{\mathbf{A}}_n + O_P\{h^{p+1} + c_n\lh\}
\end{equation}
uniformly in $u\in\Omega$. Hence by the Convexity Lemma
(\cite{Pollard(1991)}), equation (\ref{eqn:maxeqn}) also holds
uniformly in $\bstar \in C$ for any compact set $C$. Using Lemma
A.1 of \cite{Carrolletel(1997)}, it yields that
\begin{equation}\label{eqn:uniform}
    \sup_{u \in \Omega}|\bhatstar - \tilde{\mathbf{A}}_n^{-1}\mathbf{W}_n|
    \pcv 0.
\end{equation}
Furthermore, by its definition, $\bhatstar$ solves the local
likelihood equation:
\begin{equation}\label{eqn:bstarsolves}
   \sum_{i=1}^n q_1(\abar + c_n\X_i^{*T}\bhatstar,
    Y_{ni})\X_i^*K_h(U_i-u) = 0.
\end{equation}
Expanding $q_1(\abar + c_n\X_i^{*T}\bhatstar,\cdot)$ at $\abar$
yields
\begin{equation}\label{eqn:expanded}
    \mathbf{W}_n + \mathbf{A}_n\bhatstar +
    \frac{hc_n^3}{2}\sum_{i=1}^n q_3(\abar + \hat{\zeta}_i,
    Y_{ni})\X_i^*(\X_i^{*T}\bhatstar)^2K_h(U_i-u) = 0
\end{equation}
where $\hat{\zeta}_i$ lies between 0 and $c_n\X_i^{*T}\bhatstar$.
Using Conditions (A) and (C), the last term has order
$O_P(c_n^3hn\|\bhatstar\|^2)=O_P(c_n\|\bhatstar\|^2)$. 
With this, combining (\ref{eqn:An}) and (\ref{eqn:expanded}), we
obtain
\begin{equation}\label{eqn:bhatstar}
\bhatstar = \tilde{\mathbf{A}}_n^{-1}\mathbf{W}_n + O_P(h^{p+1} +
c_n\log^{1/2}(1/h))
\end{equation}
holds uniformly in $u\in\Omega$ by (\ref{eqn:uniform}). Using the
result of \cite{MackSilverman(1982)} on $\mathbf{W}_n$, we obtain
\begin{equation}\label{eqn:order}
    \|\hat{\mathbf{a}}_{0\b}(u) - \a_\b(u) \| = O_P(h^{p+1} +
c_n\lh)
\end{equation}
which holds uniformly in $u\in\Omega$.

Differentiate both sides of (\ref{eqn:bstarsolves}) w.r.t.
$\beta_{nj}$,
\begin{equation}\label{eqn:A1}
\sum_{i=1}^n q_2(\abar + c_n\X_i^{*T}\bhatstar, Y_{ni}) \biggl \{
\frac{\partial\abar}{\partial\beta_{nj}} +
c_n\biggl(\frac{\partial\bhatstar}{\partial\beta_{nj}}
\biggr)^T\X_i^*\biggr \} \X_i^*K_h(U_i-u)=0,
\end{equation}
which holds for all $u\in\Omega$. By Taylor's expansion and
similar treatments to (\ref{eqn:expanded}),
$$ \mathbf{W}_n^1 + \mathbf{W}_n^2 + (\mathbf{A}_n +
\mathbf{B}_n^1 +
\mathbf{B}_n^2)\frac{\partial\bhatstar}{\partial\beta_{nj}} +
O_P(c_n\|\bhatstar\|^2)=0,
$$ where
\begin{eqnarray*}
  \mathbf{W}_n^1 &=& hc_n\sum_{i=1}^n q_2(\abar, Y_{ni})
\frac{\partial\abar}{\partial\beta_{nj}} \X_i^*K_h(U_i-u), \\
  \mathbf{W}_n^2 &=& hc_n\sum_{i=1}^n q_3(\abar, Y_{ni})c_n\X_i^{*T}\bhatstar
\frac{\partial\abar}{\partial\beta_{nj}} \X_i^*K_h(U_i-u), \\
  \mathbf{B}_n^1 &=& hc_n^2\sum_{i=1}^n q_3(\abar, Y_{ni})c_n\X_i^{*T}\bhatstar
 \X_i^*\X_i^{*T}K_h(U_i-u), \\
  \mathbf{B}_n^2 &=& \frac{hc_n^2}{2}\sum_{i=1}^n q_4(\abar + \hat{\zeta}_i,
  Y_{ni})(c_n^2\X_i^{*T}\bhatstar)^2
 \X_i^*\X_i^{*T}K_h(U_i-u),
\end{eqnarray*}
with $\hat{\zeta}_i $ lies between 0 and $c_n\X_i^{*T}\bhatstar $.
The above equations hold for all $u\in\Omega$.
The order of $\mathbf{W}_n^2$ is smaller than that of
$\mathbf{W}_n^1$, and the orders of $\mathbf{B}_n^1$ and
$\mathbf{B}_n^2$ are smaller than that of $\mathbf{A}_n$. Hence
\begin{equation*}
\frac{\partial\bhatstar}{\partial\beta_{nj}} =
\tilde{\mathbf{A}}_n^{-1}\mathbf{W}_n^1 + o_P(\lh +
c_n^{-1}h^{p+1})
\end{equation*}
uniformly in $u\in\Omega$. From this, for $j=1,\cdots,p_n$, we
have
\begin{equation}\label{eqn:dorder}
\biggl\|\frac{\partial\hat{\mathbf{a}}_{0\b}(u)}{\partial\beta_{nj}}
- \frac{\partial\a_\b(u)}{\partial\beta_{nj}} \biggr\| =
O_P(h^{p+1} + c_n\lh)
\end{equation}
uniformly in $u\in\Omega$. Differentiating (\ref{eqn:bstarsolves})
again w.r.t. $\beta_{nk}$ and repeating as needed, we get the
desired results for higher order derivatives by following similar
arguments as above. $\square$

\begin{lemma}\label{lemma:order4}
Under Regularity Conditions (A) - (G), if $p_n^s/n \rightarrow 0$
for $s > 5/4$, $ h=O(n^{-a}) $ with $ (2s(p+1))^{-1} < a <
1-s^{-1} $, then for each $\b\in\Omega_n$,
$$ n^{-1/2}\|\nabla\hat{Q}_n(\b) - \nabla Q_n(\b)\| =
o_P(1).
$$
\end{lemma}

\noindent {\bf Proof of Lemma \ref{lemma:order4}.} Define
\begin{equation*}
\begin{split}
\mathbf{K}_1 &= n^{-1/2}\sum_{i=1}^n q_2({m}_{ni}(\b),
Y_{ni})(\Z_{ni} + {\a}_{\b}^{\prime}(U_i)\X_i)(\ahat(U_i) -
\a_\b(U_i))^T\X_i,\\
\mathbf{K}_2 &= n^{-1/2}\sum_{i=1}^n q_1({m}_{ni}(\b),
Y_{ni})(\ahatp(U_i) - \ap(U_i))^T \X_i.
\end{split}
\end{equation*}
Then by Taylor's expansion, Lemma \ref{lemma:order1} and Condition
(C),
$$
n^{-1/2} (\nabla\hat{Q}_n(\b) -  \nabla Q_n(\b)) = \mathbf{K}_1 +
\mathbf{K}_2 + \text{smaller order terms}.
$$
Define, for $\Omega$ as in Condition (F),
$$ S=\{f\in C^2(\Omega):\|f\|_{\infty} \leq 1 \}, $$
equipped with a metric $\rho(f_1,f_2) = \|f_1 - f_2\|_{\infty}$,
where $\|f\|_{\infty} = \sup_{u\in\Omega}|f(u)|$. We also let, for
$r=1,\cdots,q$ and $l=1,\cdots,p_n$,
\begin{eqnarray*}
  A_{rl}(y,u, \mathbf{X}, \mathbf{Z}_n) &=& q_2(\mathbf{X}^T\a_{\b}(u) + \mathbf{Z}_n^T\b,y)X_r\biggl(Z_{nl} + \mathbf{X}^T\frac{\partial\a_\b(u)}{\partial\beta_{nl}}\biggr), \\
  B_{r}(y,u, \mathbf{X}, \mathbf{Z}_n) &=& q_1(\mathbf{X}^T\a_{\b}(u) + \mathbf{Z}_n^T\b,y)X_r.
\end{eqnarray*}

By Lemma \ref{lemma:order1}, for any positive sequences
$(\delta_n)$ with $\delta_n \rightarrow 0$ as $n \rightarrow
\infty$, we have $P_0(\lambda_r \in S) \rightarrow 1$ and
\mbox{$P_0(\gamma_{rl} \in S) \rightarrow 1$}, where
\begin{equation*}
\begin{split}
\lambda_r &= \delta_n(h^{p+1} + c_n\lh
)^{-1}(\hat{\alpha}_{\b}^{(r)} -
\alpha_{\b}^{(r)}),\\
\gamma_{rl} &= \delta_n(h^{p+1} + c_n\lh )^{-1}
\biggl(\frac{\partial\hat{\alpha}_{\b}^{(r)}}{\partial\beta_{nl}}
- \frac{\partial\alpha_{\b}^{(r)}}{\partial\beta_{nl}}\biggr),
\end{split}
\end{equation*}
$r=1,\cdots,q$ and $l=1,\cdots,p_n.$ Hence for sufficiently large
$n$, we have $\lambda_r, \gamma_{rl}\in S$. The following three
points allow us to utilize \cite{JainMarcus(1975)} to prove our
lemma.

\begin{itemize}
    \item [I.] For any $v \in S$, we will view the map $v \mapsto A_{rl}(y,u, \mathbf{X}, \mathbf{Z}_n)v(u)$ as an element of
$C(S)$, the space of continuous functions on $S$ equipped with the
sup norm. For $v_1, v_2 \in S$, we have
\begin{equation*}
\begin{split}
&|A_{rl}(y,u, \mathbf{X}, \mathbf{Z}_n)v_1(u) - A_{rl}(y,u,
\mathbf{X}, \mathbf{Z}_n)v_2(u)|\\ &= |A_{rl}(y,u, \mathbf{X},
\mathbf{Z}_n)(v_1-v_2)(u)| \leq |A_{rl}(y,u, \mathbf{X},
\mathbf{Z}_n)|\|v_1 - v_2\|.
\end{split}
\end{equation*}
Similar result holds for $B_{r}(y,u, \mathbf{X}, \mathbf{Z}_n)$.

    \item [II.] Note that equation (\ref{eqn:defining}) is true for all
$\b$, and by differentiating w.r.t. $\b$ we get the following
formulas:
\begin{equation*}
\begin{split}
&E_0(q_1(m_{n}(\b), Y_{n})\X|U=u) = \mathbf{0},\\
&E_0(q_2(m_{n}(\b), Y_{n})\X(\Z_{n} + \ap(U)\X)^T|U=u)=\mathbf{0}.
\end{split}
\end{equation*}
Thus, we can easily see that
$$E_0(A_{rl}(Y,U, \mathbf{X}, \mathbf{Z}_n))=0$$ for each $r=1,\cdots,q$ and $l=1,\cdots,p_n$. Also we have
$$E_0(A_{rl}(Y,U, \mathbf{X}, \mathbf{Z}_n)^2) < \infty, $$ by Regularity Conditions (A) and (C).
For $B_{r}(Y,U, \mathbf{X}, \mathbf{Z}_n)$, results hold
similarly.
    \item [III.] Let $H(\cdot,S)$ denote the metric entropy of the set $S$ w.r.t. the metric $\rho$. Then
$$H(\epsilon, S) \leq C_0\epsilon^{-1} $$ for some constant $C_0$. Hence
$\int_0^1 H^{1/2}(\epsilon, S)d\epsilon < \infty.$
\end{itemize}

Conditions of Theorem 1 in \cite{JainMarcus(1975)} can be derived
from the three notes above, so that we have
$$ n^{-1/2}\sum_{i=1}^n A_{rl}(Y_i,U_i, \mathbf{X}_i, \mathbf{Z}_{ni})(\cdot),$$
where $A_{rl}(Y_i,U_i, \mathbf{X}_i, \mathbf{Z}_{ni})(\cdot), \;
i=1,\cdots,n$ being i.i.d. replicates of $A_{rl}(Y,U, \mathbf{X},
\mathbf{Z}_n)(\cdot)$ in $C(S)$, converges weakly to a Gaussian
measure on $C(S)$. Hence, since $\lambda_r, \gamma_{rl}\in S$,
$$ n^{-1/2}\sum_{i=1}^n A_{rl}(Y_i,U_i, \mathbf{X}_i, \mathbf{Z}_{ni})(\lambda_r) = O_P(1), $$ which
implies that
$$ n^{-1/2}\sum_{i=1}^n A_{rl}(Y_i,U_i, \mathbf{X}_i, \mathbf{Z}_{ni})(\hat{\alpha}_{\b}^{(r)} -
\alpha_{\b}^{(r)}) = O_P(\delta_n^{-1}(h^{p+1} + c_n\lh ) ).$$
Similarly, apply Theorem 1 of \cite{JainMarcus(1975)} again, we
have
$$ n^{-1/2}\sum_{i=1}^n B_{r}(Y_i,U_i, \mathbf{X}_i, \mathbf{Z}_{ni})
\biggl (\frac{\partial\hat{\alpha}_{\b}^{(r)}}{\partial\beta_{nl}}
- \frac{\partial\alpha_{\b}^{(r)}}{\partial\beta_{nl}}\biggr )
=O_P(\delta_n^{-1}(h^{p+1} + c_n\lh ) ).$$ Then the column vector
$\mathbf{K}_1$ which is $p_n-$dimensional, has the $l^{\text{th}}$
component equals
$$\sum_{r=1}^q\biggl \{n^{-1/2}\sum_{i=1}^n A_{rl}(Y_i,U_i, \mathbf{X}_i, \mathbf{Z}_{ni})(\hat{\alpha}_{\b}^{(r)} -
\alpha_{\b}^{(r)})\biggr \} = O_P(\delta_n^{-1}(h^{p+1} + c_n\lh
),
$$ using the result just proved. Hence we have shown
$$\|\mathbf{K}_1\| = O_P(\sqrt{p_n}\delta_n^{-1}(h^{p+1} +
 c_n\lh ))  = o_P(1),$$
since $\delta_n$ can be made arbitrarily slow in converging to 0.
Similarly, we have $\|\mathbf{K}_2\| = o_P(1)$ as well. The
conclusion of the lemma follows.  $\;\square$

\vspace{12pt}
\noindent {\bf Proof of Theorem \ref{thm:ple}}.

Let $\gamma_n = \sqrt{p_n/n}$. Our aim is to show that, for a
given $\epsilon > 0$,
\begin{equation}\label{eqn:consistency}
\mathbb{P}\biggl\{\sup_{\|\mathbf{v}\|=C}\hat{Q}_n(\b_0 +
\gamma_n\mathbf{v}) < \hat{Q}_n(\b_0) \biggr\}\geq 1-\epsilon,
\end{equation}
so that this implies with probability tending to 1 there is a
local maximum $\bhat$ in the ball $\{\b_0 + \gamma_n\mathbf{v} :
\|\mathbf{v}\| \leq C\} $ such that $\|\bhat - \b_0 \| =
O_P(\gamma_n)$.

Define the terms $\hat{I}_1 = \g \nabla^T\Q(\b_0)\v, \hat{I}_2 =
\frac{\gsq}{2}\v^T\nabla^2\Q(\b_0)\v$ and \\
$\hat{I}_3 =
\frac{\gcu}{6}\nabla^T(\v^T\nabla^2\Q(\boldsymbol{\beta}_n^*)\v)\v$.
By Taylor's expansion,
\begin{eqnarray*}
 \hat{Q}_n(\b_0 + \gamma_n\mathbf{v}) -
\hat{Q}_n(\b_0) = \hat{I}_1 + \hat{I}_2 + \hat{I}_3,
\end{eqnarray*}
where $\boldsymbol{\beta}_n^*$ lies between $\b_0$ and $\b_0 +
\g\v$.

We further split $\hat{I}_1 = D_1 + D_2$, where
\begin{equation*}
\begin{split}
           D_1 &= \sum_{i=1}^n q_1(\hat{m}_{ni}(\b_0), Y_{ni})(\Z_{ni} + \aop(U_i)\X_i)^T\v\g, \\
           D_2 &= \sum_{i=1}^n q_1(\hat{m}_{ni}(\b_0),Y_{ni})\X_i^T(\ahatop(U_i) - \aop(U_i))^T\v\g,
\end{split}
\end{equation*}
with $\hat{m}_{ni}(\b) = \ahat(U_i)^T\X_i +
\boldsymbol{\beta}_n^T\Z_{ni}$. By Condition (A) and Lemma
\ref{lemma:order1}, $D_2$ has order smaller than $D_1$. Using
Taylor's expansion, we have
$$D_1 = \g\v^T \biggl ( \sum_{i=1}^n \frac{\partial Q_{ni}(\b_0)}
{\partial\b} + \sqrt{n}\mathbf{K}_1 \biggr ) + \text{smaller order
terms},$$ where $\mathbf{K}_1$ is as defined in Lemma
\ref{lemma:order4} so that within the lemma's proof we have
$\|\mathbf{K}_1\| = o_P(1)$. Using equation (\ref{eqn:bartlett}),
we have by the mean-variance decomposition
$$ \biggl\|\v^T\sum_{i=1}^n \frac{\partial Q_{ni}(\b_0)}{\partial\b}\biggr\|
= O_P(\sqrt{n\v^TI_n(\b_0)\v}) = O_P(\sqrt{n}) \|\v\|, $$ where
last inequality follows from Condition (B).
Hence
\begin{equation*}
|\hat{I}_1| = O_P(\sqrt{n}\gamma_n) \|\v\|.
\end{equation*}

Next, consider $\hat{I}_2 = I_2 + (\hat{I}_2 - I_2)$, where
\begin{equation*}
\begin{split}
I_2 &= \frac{1}{2}\v^T\nabla^2Q_n(\b_0)\v\gsq\\
    &= - \frac{n}{2}\v^TI_n(\b_0)\v\gsq +
    \frac{n}{2}\v^T\{n^{-1}\nabla^2Q_n(\b_0) + I_n(\b_0)
    \}\v\gsq \\
    &=- \frac{n}{2}\v^TI_n(\b_0)\v\gsq +
    o_P(n\gsq)\|\v\|^2
\end{split}
\end{equation*}
with the last line follows from Lemma \ref{lemma:1} in the
Appendix. Using Lemma \ref{lemma:8}, $$\| \hat{I}_2 - I_2 \| =
o_P(n\gsq\|\v\|^2).$$

On the other hand, by Condition (B), we have
\begin{equation*}
    |n\gsq\v^TI_n(\b_0)\v | \geq
    O(n\gsq\lambda_{\text{min}}(I_n(\b_0))\|\v\|^2)=O(n\gsq\|\v\|^2).
\end{equation*}
Hence, $\hat{I}_2 - I_2$ has a smaller order than $I_2$.

Finally consider $\hat{I}_3$. We suppress the dependence of
$\a_\b(U_i)$ and  its derivatives on $U_i$, and denote $q_{1i} =
q_1(m_{ni}(\b_0), Y_{ni})$. Using Taylor's expansions, expanding
$\hat{Q}_n(\bbeta_n^*)$ at $\b_0$ and then $\hat{Q}_n(\b_0)$ at
$\a_{\b_0}$, we can arrive at
\begin{equation*}
\begin{split}
\hat{Q}_n(\bbeta_n^*) = Q_n(\b_0) &+  \sum_{i=1}^n  \{
q_{1i}\mathbf{X}_i^T (\ahato - \ao)\\ &+ q_{1i}(\mathbf{Z}_{ni} +
\ahatop\mathbf{X}_i)^T (\bbeta_n^* - \b_0)  \}(1+o_P(1)).
\end{split}
\end{equation*}
Substituting $\hat{Q}_n(\bbeta_n^*)$ into $\hat{I}_3$ with the
right hand side above, by Condition (C) and Lemma
\ref{lemma:order1}, we have
\begin{equation*}
\begin{split}
\hat{I}_3 &=
\frac{1}{6}\sum_{i,j,k=1}^{p_n}\frac{\partial^3Q_n(\boldsymbol{\beta}_{n0})}{\partial\b_i\partial\b_j\partial\b_k}
v_iv_jv_k\gcu + \text{smaller order terms.}\\
\end{split}
\end{equation*}
Hence,
\begin{equation*}
|\hat{I}_3| = O_P(np_n^{3/2}\gcu\|\v\|^3)
=O_P(\sqrt{p_n^4/n}\|\v\|)n\gsq\|\v\|^2 = o_P(1)n\gsq\|\v\|^2.
\end{equation*}
Comparing, we find the order of $-n\gsq \v^TI_n(\b_0)\v$ dominates
all other terms by allowing $\|\v\|=C$ to be large enough. This
proves (\ref{eqn:consistency}). $\square$

\vspace{12pt}
\noindent {\bf Proof of Theorem \ref{thm:asymnorm}.}

Note that by Theorem \ref{thm:ple}, $\|\bhat - \b_0\| =
O_P(\sqrt{p_n/n})$. Since $\nabla
\hat{Q}_n(\hat{\boldsymbol{\beta}}_n)=0$, by Taylor's expansion,
\begin{equation} \label{eqn:simplify}
\nabla \hat{Q}_n(\b_0) + \nabla^2
\hat{Q}_n(\b_0)(\hat{\boldsymbol{\beta}}_n - \b_0) + \mathcal{C} =
0,
\end{equation}
where $\boldsymbol{\beta}_n^*$ lies between $\b_0$
$\hat{\boldsymbol{\beta}}_n$  and $\mathcal{C} =
\frac{1}{2}(\hat{\boldsymbol{\beta}_n} - \b_0)^T\nabla^2(\nabla
\hat{Q}_n(\boldsymbol{\beta}_n^*))(\hat{\boldsymbol{\beta}}_n -
\b_0))$ which is understood as a vector of quadratic components.

Using similar argument to approximating $\hat{I}_3$ in Theorem
\ref{thm:ple}, by Lemma \ref{lemma:order1} and noting
$\|\boldsymbol{\beta}_n^* - \b_0\| = o_P(1)$, we have $\|\nabla^2
\frac{\partial
\hat{Q}_n(\boldsymbol{\beta}_n^*)}{\partial\beta_{nj}} \|^2 =
O_P(n^2p_n^2)$. Hence
\begin{equation}\label{eqn:heavyassumption}
\|n^{-1}\mathcal{C}\|^2 \leq n^{-2} \|\hat {\boldsymbol{\beta}}_n
- \b_0\|^4 \sum_{j=1}^{p_n} \biggl \|\nabla^2 \frac{\partial
\hat{Q}_n(\bbeta_n^*)}{\partial\beta_{nj}} \biggr \|^2
=O_P(p_n^5/n^2) = o_P(n^{-1}).
\end{equation}
At the same time, by Lemma \ref{lemma:1} and the Cauchy-Schwarz
inequality,
\begin{equation}\label{eqn:simplify2}
\begin{split}
&\|n^{-1}\nabla^2\hat{Q}_n(\b_0) (\bhat - \b_0) +
I_n(\b_0)(\bhat - \b_0) \|\\
&= o_P((np_n)^{-1/2}) + O_P(\sqrt{p_n^3/n}(h^{p+1} + c_n\lh)) =
o_P(n^{-1/2}).
\end{split}
\end{equation}
Combining (\ref{eqn:simplify}),(\ref{eqn:heavyassumption}) and
(\ref{eqn:simplify2}), we have
\begin{equation}\label{eqn:approx}
\begin{split}
I_n(\b_0)(\bhat - \b_0) &= n^{-1}\nabla \hat{Q}_n(\b_0) +
o_P(n^{-1/2})\\
&= n^{-1}\nabla Q_n(\b_0) + o_P(n^{-1/2}),
\end{split}
\end{equation}
where the last line follows from Lemma \ref{lemma:order4}.
Consequently, using equation (\ref{eqn:approx}), we get
\begin{equation}\label{eqn:norm1}
\begin{split}
\sqrt{n} A_nI_n^{1/2}(\b_0)(\bhat - \b_0)
&=n^{-1/2}A_nI_n^{-1/2}(\b_0)\nabla Q_n(\b_0) \\&\;\;\;\;+
o_P(A_nI_n^{-1/2}(\b_0))\\
&=n^{-1/2}A_nI_n^{-1/2}(\b_0)\nabla Q_n(\b_0) + o_P(1),
\end{split}
\end{equation}
since  $\|A_nI_n^{-1/2}(\b_0)\|$ = $O(1)$ by conditions of Theorem
\ref{thm:asymnorm}.

We now check the Lindeberg-Feller Central Limit Theorem (see for
example, \cite{vanderVaart(1998)}) for the last term in
(\ref{eqn:norm1}). Let $B_{ni} = n^{-1/2}A_nI_n^{-1/2}(\b_0)\nabla
Q_{ni}(\b_0)$, \mbox{$i=1,\cdots,n$}. Given $\epsilon>0$,
\begin{equation*}
\sum_{i=1}^n E_0\|B_{ni}\|^21\{\|B_{ni}\| > \epsilon \}
\leq n\sqrt{E_0\|B_{n1}\|^4 \cdot \mathbb{P}(\|B_{n1}\| >
\epsilon)}.
\end{equation*}
Using Chebyshev's inequality,
\begin{equation}\label{eqn:estforprob}
\begin{split}
\mathbb{P}(\|B_{n1}\|>\epsilon) &\leq
n^{-1}\epsilon^{-2}E\|A_nI_n^{-1/2}(\b_0)\nabla
Q_{n1}(\b_0)\|^2\\
&=n^{-1}\epsilon^{-2}tr(G) = O(n^{-1}),
\end{split}
\end{equation}

\noindent where $tr(A)$ is the trace of square matrix A.
Similarly, we can show that, using Condition~(B),
\begin{equation}\label{eqn:estforexp}
\begin{split}
E_0\|B_{n1}\|^4 &\leq
\sqrt{l}n^{-2}\lambda_{\text{min}}^2(A_nA_n^T)\lambda_{\text{max}}^2(I_n(\b_0))\sqrt{E_0(\nabla
Q_{n1}(\b_0)^T \nabla Q_{n1}(\b_0))^4}\\
&=O(p_n^2/n^2).
\end{split}
\end{equation}
Therefore (\ref{eqn:estforprob}) and (\ref{eqn:estforexp})
together imply
\begin{equation*}
\sum_{i=1}^nE_0\|B_{ni}\|^21\{\|B_{ni}\| > \epsilon \}  = O(
\sqrt{p_n^2/n} ) = o(1).
\end{equation*}
Also,
\begin{equation*}
\begin{split}
\sum_{i=1}^n \text{Var}_0(B_{ni}) & =
\text{Var}_0(A_nI_n^{-1/2}(\b_0)\nabla Q_{n1}(\b_0))\\ &=A_nA_n^T
\rightarrow G.
\end{split}
\end{equation*}
Therefore $B_{ni}$ satisfies the conditions of the
Lindeberg-Feller Central Limit Theorem. Consequently, using
(\ref{eqn:norm1}), it follows that
$$ \sqrt{n}A_nI_n^{1/2}(\b_0)(\bhat - \b_0) \dcv N(0, G),$$
and this completes the proof. $\thickspace \square$


Referring back to Section \ref{subsect:hyptest}, let $B_n$ be a
$(p_n - l) \times p_n$ matrix satisfying $B_nB_n^T = I_{p_n - l}$
and $A_nB_n^T = 0$. Since $A_n\b=0$ under $H_0$, rows of $A_n$ are
perpendicular to $\b$ and the orthogonal complement of rows of
$A_n$ is spanned by rows of $B_n$ since $A_nB_n^T=0$. Hence
$$\b = B_n^T\boldsymbol{\gamma}$$ under $H_0$, where
$\boldsymbol{\gamma}$ is a $(p_n - l) \times 1$ vector. Then under
$H_0$ the profile likelihood estimator is also the local maximizer
$\boldsymbol{\hat{\gamma}}_n$ of the problem
$$ \hat{Q}_n (B_n^T\ghatvec) = \max_{\gvec}Q_n(B_n^T\gvec).$$

\vspace{12pt} \noindent {\bf Proof of Theorem \ref{thm:hyptest}.}

By Taylor's expansion, expanding $\hat{Q}_n(B_n^T\ghatvec)$ at
$\bhat$ and noting that $\nabla^T \hat{Q}_n(\bhat)=0$, then
$\hat{Q}_n(\bhat) - \hat{Q}_n(B_n^T\ghatvec)=T_1 + T_2,$ where
\begin{equation*}
\begin{split}
T_1 &= -\frac{1}{2}(\bhat - B_n^T\ghatvec)^T\nabla^2\hat{Q}_n(\bhat)(\bhat - B_n^T\ghatvec),\\
T_2 &= \frac{1}{6}\nabla^T\{(\bhat -
B_n^T\ghatvec)^T\nabla^2\hat{Q}_n(\beta_n^*)(\bhat -
B_n^T\ghatvec)\}(\bhat - B_n^T\ghatvec).
\end{split}
\end{equation*}
Denote by $\Theta_n = I_n(\b_0)$ and $\boldsymbol{\Phi}_n =
\frac{1}{n}\nabla Q_n(\b_0)$. Using equation (\ref{eqn:approx})
and noting that $\Theta_n$ has eigenvalues uniformly bounded away
from 0 and infinity by Condition (B), we have
$$\bhat - \b_0 = \Theta_n^{-1}\boldsymbol{\Phi}_n  +
o_P(n^{-1/2}).$$ Combining this with Lemma \ref{lemma:4} in the
Appendix, under the null hypothesis $H_0$,
\begin{equation}\label{eqn:lemma5eqn1}
\begin{split}
\bhat - B_n^T\ghatvec = &\Theta_n^{-1/2}\{I_{p_n} - \Theta_n^{1/2}
B_n^T(B_n\Theta_nB_n^T)^{-1}B_n\Theta_n^{1/2}\}
\Theta_n^{-1/2}\boldsymbol{\Phi}_n \\
&+ o_P(n^{-1/2}).
\end{split}
\end{equation}

Since $S_n = I_{p_n} - \Theta_n^{1/2}B_n^T
(B_n\Theta_nB_n^T)^{-1}B_n\Theta_n^{1/2}$ is a $p_n \times p_n$
idempotent matrix with rank $l$, it follows by mean-variance
decomposition of the term $\| \bhat - B_n^T\ghatvec \|^2$ and
Condition (B) that
$$\| \bhat - B_n^T\ghatvec \| = O_P(n^{-1/2}).$$
Hence, using similar argument as in the approximation of order for
$|\hat{I}_3|$ in Theorem \ref{thm:ple}, we have
\begin{equation*}
\begin{split}
|T_2| = O_P(np_n^{3/2}) \cdot \|\bhat - B_n^T\ghatvec\|^3 =o_P(1).
\end{split}
\end{equation*}
Hence $\hat{Q}_n(\bhat) - \hat{Q}(B_n^T\ghatvec) = T_2 + o_P(1)$.

By Lemma \ref{lemma:1} and an approximation to $n^{-1}\|\nabla^2
\hat{Q}_n(\bhat) - \nabla^2 \hat{Q}_n(\b_0)\| = o_P( p_n^{-1/2} )$
(the proof is similar to that for Lemma \ref{lemma:7} with the
proof of order for $|\hat{I}_3|$ in Theorem \ref{thm:ple}, and is
omitted), we have
\begin{equation*}
\begin{split}
&\biggl\|\frac{1}{2}(\bhat - B_n^T\ghatvec)^T\{\nabla^2\hat{Q}_n(\bhat) + nI_n(\b_0)\}(\bhat - B_n^T\ghatvec)\biggr\| \\
&= O_P( l/n )\cdot n\{ o_P( p_n^{-1/2} ) + O_P( p_n(h^{p+1} +
c_n\lh) )\} = o_p(1).
\end{split}
\end{equation*}
Therefore,
$$
   \hat{Q}_n(\bhat) - \hat{Q}_n(B_n^T\ghatvec)
   = \frac{n}{2}(\bhat - B_n^T\ghatvec)^TI_n(\b_0)(\bhat - B_n^T\ghatvec) + o_P(1).
$$
By (\ref{eqn:lemma5eqn1}), we have
$$
 \hat{Q}_n(\bhat) - \hat{Q}_n(B_n^T\ghatvec) =
 \frac{n}{2}\boldsymbol{\Phi}_n^T\Theta_n^{-1/2}
 S_n\Theta_n^{-1/2}\boldsymbol{\Phi}_n + o_P(1).
$$

Since $S_n$ is idempotent, it can be written as $S_n = D_n^TD_n$
where $D_n$ is an $l \times p_n$ matrix satisfying $D_nD_n^T =
I_l$. By Theorem \ref{thm:asymnorm},  we have already shown that
$\sqrt{n}D_n\Theta_n^{-1/2}\boldsymbol{\Phi}_n \dcv N(\mathbf{0},
I_l)$. Hence
$$ 2\{\hat{Q}_n(\bhat) - \hat{Q}_n(B_n^T\ghatvec)\} =
n(D_n\Theta_n^{-1/2}\boldsymbol{\Phi}_n)^T(D_n\Theta_n^{-1/2}\boldsymbol{\Phi}_n)
\dcv \chi_l^2. \thickspace \square $$

\vspace{12pt}
\noindent {\bf Proof of Theorem \ref{thm:sandwich}.}

Let $\hat{\mathcal{A}}_n = -n^{-1}\nabla^2 \hat{Q}_n(\bhat)$,
$\hat{\mathcal{B}}_n =
\widehat{\text{cov}}\{\nabla\hat{Q}_n(\bhat)\}$  and  $\mathcal{C}
= I_n(\b_0)$. Write
$$I_1 = \hat{\mathcal{A}}_n^{-1}(\hat{\mathcal{B}}_n - \mathcal{C})\hat{\mathcal{A}}_n^{-1}, \thickspace\thickspace\thickspace\thickspace
I_2 = \hat{\mathcal{A}}_n^{-1}(\mathcal{C} -
\hat{\mathcal{A}}_n)\hat{\mathcal{A}}_n^{-1},\thickspace\thickspace\thickspace\thickspace
I_3 = \hat{\mathcal{A}}_n^{-1}(\mathcal{C} - \hat{\mathcal{A}}_n)
\mathcal{C}^{-1}.
$$
Then, $ \hat{\Sigma}_n - \Sigma_n = I_1 + I_2 + I_3$. Our aim is
to show that, for all $i = 1,\cdots,p_n$,
$$\lambda_i(\hat{\Sigma}_n - \Sigma_n) = o_P(1),$$ so that
$A_n(\hat{\Sigma}_n - \Sigma_n)A_n^T \pcv 0$, where $\lambda_i(A)$
is the $i$th eigenvalue of a symmetric matrix A. Using the
inequalities
\begin{equation*}
\begin{split}
\lambda_{\text{min}}(I_1)+\lambda_{\text{min}}(I_2)+\lambda_{\text{min}}(I_3)
&\leq
\lambda_{\text{min}}(I_1+I_2+I_3)\\
 \lambda_{\text{max}}(I_1 + I_2 + I_3) & \leq
\lambda_{\text{max}}(I_1)+\lambda_{\text{max}}(I_2)+\lambda_{\text{max}}(I_3),\\
\end{split}
\end{equation*}
it suffices to show that $\lambda_i(I_j)=o_P(1)$ for $j=1,2,3.$
From the definition of $I_1, I_2$ and $I_3$, it is clear that we
only need to show $\lambda_i(\mathcal{C}-\hat{\mathcal{A}}_n) =
o_P(1)$ and $\lambda_i(\hat{\mathcal{B}}_n - \mathcal{C}) =
o_P(1)$. Let $K_1 = I_n(\b_0) + n^{-1}\nabla^2Q_n(\b_0)$, $K_2 =
n^{-1}(\nabla^2Q_n(\bhat) - \nabla^2Q_n(\b_0))$, and $K_3 =
n^{-1}(\nabla^2\hat{Q}_n(\bhat) - \nabla^2Q_n(\bhat))$. Then,
 $$
    \mathcal{C} - \hat{\mathcal{A}}_n = K_1 + K_2 + K_3.
 $$
Applying Lemma \ref{lemma:1} to $K_1$, Lemma \ref{lemma:7} to
$K_2$, and Lemma \ref{lemma:8} to $K_3$, we have $\|\mathcal{C} -
\hat{\mathcal{A}}\| = o_P(1)$. Thus, $\lambda_i(\mathcal{C} -
\hat{\mathcal{A}}) = o_P(1)$. Hence the only thing left to show is
$\lambda_i(\hat{\mathcal{B}}_n - \mathcal{C}) = o_P(1)$.

To this end, consider the decomposition
$$\hat{\mathcal{B}}_n - \mathcal{C} = K_4 + K_5$$
where
\begin{eqnarray*}
  K_4 &=& \biggl\{\frac{1}{n}\sum_{i=1}^n \frac{\partial\hat{Q}_{ni}(\bhat)}{\partial\beta_{nj}}
  \frac{\partial\hat{Q}_{ni}(\bhat)}{\partial\beta_{nk}} \biggr\} - I_n(\b_0), \\
  K_5 &=& -  \biggl\{\frac{1}{n}\sum_{i=1}^n \frac{\partial\hat{Q}_{ni}(\bhat)}{\partial\beta_{nj}}\biggr\}
         \biggl\{\frac{1}{n}\sum_{i=1}^n \frac{\partial\hat{Q}_{ni}(\bhat)}{\partial\beta_{nk}}\biggr\}.
\end{eqnarray*}
Our goal is to show that $K_4$ and $K_5$ are $o_P(1)$, which then
implies $\lambda_i(\hat{\mathcal{B}}_n - \mathcal{C}) = o_P(1)$.
We consider $K_4$ first, which can be further decomposed into  $
K_4 = K_6 + K_7$, where
\begin{eqnarray*}
  K_6 &=& \biggl\{\frac{1}{n}\sum_{i=1}^n \frac{\partial\hat{Q}_{ni}(\bhat)}{\partial\beta_{nj}}
\frac{\partial\hat{Q}_{ni}(\bhat)}{\partial\beta_{nk}} -
\frac{1}{n}\sum_{i=1}^n \frac{\partial
Q_{ni}(\b_0)}{\partial\beta_{nj}}
\frac{\partial Q_{ni}(\b_0)}{\partial\beta_{nk}}\biggr\}, \\
  K_7 &=& \biggl\{\frac{1}{n}\sum_{i=1}^n \frac{\partial Q_{ni}(\b_0)}{\partial\beta_{nj}}
\frac{\partial Q_{ni}(\b_0)}{\partial\beta_{nk}}\biggr\} -
I_n(\b_0).
\end{eqnarray*}
Observe that
\begin{equation*}
\begin{split}
K_6 &= \bigg\{\frac{1}{n}\sum_{i=1}^n \frac{\partial
Q_{ni}(\b_0)}{\partial\beta_{nj}}
\bigg\{\frac{\partial\hat{Q}_{ni}(\bhat)}{\partial\beta_{nk}} - \frac{\partial Q_{ni}(\b_0)}{\partial\beta_{nk}}\bigg\}\\
&\thickspace\thickspace\thickspace +\frac{1}{n}\sum_{i=1}^n
\frac{\partial Q_{ni}(\b_0)}{\partial\beta_{nk}}
\bigg\{\frac{\partial\hat{Q}_{ni}(\bhat)}{\partial\beta_{nj}} - \frac{\partial Q_{ni}(\b_0)}{\partial\beta_{nj}}\bigg\}\\
&\thickspace\thickspace\thickspace +
\frac{1}{n}\sum_{i=1}^n\bigg\{\frac{\partial\hat{Q}_{ni}(\bhat)}{\partial\beta_{nk}}
- \frac{\partial
Q_{ni}(\b_0)}{\partial\beta_{nk}}\bigg\}\bigg\{\frac{\partial\hat{Q}_{ni}(\bhat)}{\partial\beta_{nj}}
- \frac{\partial Q_{ni}(\b_0)}{\partial\beta_{nj}}\bigg\}\bigg\},
\end{split}
\end{equation*}
and this suggests that an approximation of the order of
$\frac{\partial}{\partial\beta_{nk}}
(\hat{Q}_{ni}(\bhat)-Q_{ni}(\b_0))$ for each $k=1,\cdots,p_n$ and
$i=1,\cdots,n$ is rewarding. Define
\begin{eqnarray*}
  a_{ik} = \frac{\partial}{\partial\beta_{nk}}(\hat{Q}_{ni}(\bhat)-Q_{ni}(\bhat)),
 \quad \mbox{and} \quad b_{ik} = \frac{\partial}{\partial\beta_{nk}}(Q_{ni}(\bhat)-Q_{ni}(\b_0)),
\end{eqnarray*}
then
$\frac{\partial}{\partial\beta_{nk}}(\hat{Q}_{ni}(\bhat)-Q_{ni}(\b_0))
= a_{ik} + b_{ik}$. By Taylor's expansion, suppressing dependence
of $\a_\b(U_i)$ and its derivatives on $U_i$,
\begin{equation*}
a_{ik} = \biggl \{ \frac{\partial^2
Q_{ni}(\bhat)}{\partial\beta_{nk}\partial\a_\b^T} (\hat{\a}_\bhat
- \a_\bhat) + \frac{\partial Q_{ni}(\bhat)}{\partial\a_\b^T}
\biggl( \frac{\partial\hat{\a}_\bhat}{\partial\beta_{nk}} -
\frac{\partial\a_\bhat}{\partial\beta_{nk}} \biggr) \biggr \}
(1+o_P(1)).
\end{equation*}
Using Lemma \ref{lemma:order1}, Condition $(C)$, with argument
similar to the proof of Lemma \ref{lemma:8}, we then have
$$ a_{ik} = O_P(h^{p+1} + c_n\lh). $$

Similarly, Taylor's expansion gives
\begin{equation*}
\begin{split}
b_{ik} &= \frac{\partial^2
Q_{ni}(\b_0)}{\partial\beta_{nk}\partial\b^T}(\bhat -
\b_0)(1+o_P(1)),
\end{split}
\end{equation*}
which implies that, by Theorem \ref{thm:ple} and Regularity
Condition (C),
$$|b_{ik}| = O_P(\sqrt{p_n^2/n}).$$

Using the approximations of $a_{ik}$ and $b_{ik}$ above, by
Condition (C),
\begin{equation*}
\begin{split}
&\bigg|\frac{1}{n}\sum_{i=1}^n \frac{\partial
Q_{ni}(\b_0)}{\partial\beta_{nj}}\bigg\{
\frac{\partial\hat{Q}_{ni}(\bhat)}{\partial\beta_{nk}} -
\frac{\partial Q_{ni}(\b_0)}{\partial\beta_{nk}}\bigg\}\bigg|\\
&\leq \frac{1}{n}\sum_{i=1}^n\biggl|\frac{\partial Q_{ni}(\b_0)}{\partial\beta_{nj}}\biggr|\cdot|a_{ik} + b_{ik}|\\
&= O_P(h^{p+1} + c_n\lh  + n^{-1/2} p_n ).
\end{split}
\end{equation*}
This shows that
$$ \|K_6\| = O_P( p_n( h^{p+1} + c_n\lh ) + p_n^2n^{-1/2} ) = o_P(1)$$ by the conditions of the theorem.

For $K_7$, note that
\begin{equation*}
\begin{split}
E_0K_7 & = n^{-2}(np_n^2) E_0 \biggl\{\frac{\partial
Q_{ni}(\b_0)}{\partial\beta_{nj}}\frac{\partial
Q_{ni}(\b_0)}{\partial\beta_{nk}} - E_0\biggl(\frac{\partial
Q_{ni}(\b_0)}{\partial\beta_{nj}} \frac{\partial
Q_{ni}(\b_0)}{\partial\beta_{nk}}\biggr)\biggr\}^2 \\
&= O(p_n^2/n)
\end{split}
\end{equation*}
which implies that $\|(K_7)\|=O_P(p_n^2/n)=o(1)$. Hence using $K_4
= K_6 + K_7$,
$$\|K_4\| = o_P(1) + O_P ( p_n(h^{p+1} + c_n\lh ) + \sqrt{p_n^4/n}  )=o_P(1).$$

 Finally consider $K_5$. Define $A_j = n^{-1}\sum_{i=1}^n(a_{ij} + b_{ij}) + n^{-1}\sum_{i=1}^n
\frac{\partial Q_{ni}(\b_0)}{\partial\beta_{nj}}$, where $a_{ij}$
and $b_{ij}$ are defined as before,
 we can then rewrite $K_5 = \{A_jA_k\}$. Now
\begin{equation*}
\begin{split}
|A_j| &\leq \sup_{i,j}|a_{ij} + b_{ij}| +
\biggl|\frac{1}{n}\sum_{i=1}^n
\frac{\partial Q_{ni}(\b_0)}{\partial\beta_{nj}}\biggr|\\
&= O_P(h^{p+1} + c_n\lh + n^{-1/2}p_n ) + O_P(n^{-1/2}),
\end{split}
\end{equation*}
where the last line follows from the approximations for $a_{ij}$
and $b_{ij}$, and mean-variance decomposition of the term
$n^{-1}\sum_{i=1}^n \frac{\partial
Q_{ni}(\b_0)}{\partial\beta_{nj}}$. Hence
$$\|K_5\| = O_P(p_n (h^{p+1} + c_n\lh + n^{-1/2}p_n
)^2 ) =o_P(1), $$ and this completes the proof. $\; \square$

\vspace{12pt} \noindent {\bf Proof of Theorem
\ref{thm:consistentalphaprime}}.

In expression (\ref{eqn:llf}), we set $p=0$, which effectively
assumes $\a_\b(U_i) \approx \a_\b(u)$ for $U_i$ in a neighborhood
of $u$. Using the same notation as in the proof of \mbox{Lemma
\ref{lemma:order1}}, we have $\bar{\a}_{ni}(u) =
\a_\b(u)^T\mathbf{X}_i + \mathbf{Z}_{ni}^T\b$, $\bhatstar =
c_n^{-1}(\hat{\mathbf{a}}_{0\b}(u) - \a_\b(u))$ and
$\mathbf{X}_i^* = \mathbf{X}_i$. Following the proof of Lemma
\ref{lemma:order1}, we arrive at equation (\ref{eqn:A1}), which in
this case is reduced to
$$\sum_{i=1}^n q_2(\mathbf{X}_i^T\hat{\mathbf{a}}_{0\b}(u) + \mathbf{Z}_{ni}^T\b, Y_{ni})
\biggl(Z_{nij} +
\biggl(\frac{\partial\hat{\mathbf{a}}_{0\b}(u)}{\partial\beta_{nj}}\biggr)^T\mathbf{X}_i
\biggr) \mathbf{X}_iK_h(U_i-u)=0.$$ Solving for
$\frac{\partial\hat{\mathbf{a}}_{0\b}(u)}{\partial\b}$ from the
above equation, which is true for $j=1,\cdots,p_n$, we get the
same expression as given in the lemma.

 Hence it remains to show that $\frac{\partial\hat{\mathbf{a}}_{0\b}(u)}{\partial\b}$ is a consistent estimator
of $\a_{\b}^{\prime}(u)$. However this is done by the proof of
Lemma \ref{lemma:order1} already, where equation
(\ref{eqn:dorder}) becomes
$$\biggl\|\frac{\partial\hat{\mathbf{a}}_{0\b}(u)}{\partial\b} -  \hat{\a}_{\b}^{\prime}(u)\biggr\| =
O_P(\sqrt{p_n}(h + c_n\lh ))=o_P(1) $$ and the proof completes.
\thickspace $\square$

\vspace{12pt}
\centerline{\large \bf APPENDIX: PROOFS OF LEMMAS \ref{lemma:7} - \ref{lemma:4}}

\begin{lemma}\label{lemma:7}
Assuming Conditions (A) - (G) and $p_n^4/n = o(1)$, we have
$$ n^{-1}\|\nabla^2Q_n(\bhat) - \nabla^2Q_n(\b_0)\| = o_P(1). $$
\end{lemma}

\noindent \emph{Proof of Lemma \ref{lemma:7}}. Consider
\begin{equation*}
\begin{split}
  n^{-1}\|\nabla^2Q_n(\b) - \nabla^2Q_n(\b_0)\|^2
  &= \frac{1}{n^2}\sum_{i,j=1}^{p_n}
  \biggl( \frac{\partial^2Q_n(\bhat)}{\partial\beta_{ni}\partial\beta_{nj}} -
  \frac{\partial^2Q_n(\b_0)}{\partial\beta_{ni}\partial\beta_{nj}} \biggr)^2 \\
    &= \frac{1}{n^2}\sum_{i,j=1}^{p_n}\biggl(\sum_{k=1}^{p_n}\frac{\partial^3Q_n(\bstar)}
       {\partial\beta_{ni}\partial\beta_{nj}\partial\beta_{nk}}(\hat{\beta}_{nk} - \beta_{0k})\biggr)^2\\
    &\leq \frac{1}{n^2}\sum_{i,j=1}^{p_n}\sum_{k=1}^{p_n}\biggl(\frac{\partial^3Q_n(\bstar)}
       {\partial\beta_{ni}\partial\beta_{nj}\partial\beta_{nk}}\biggr)^2\|\hat{\beta}_{nk} - \beta_{0k}\|^2,
\end{split}
\end{equation*}
where $\bstar$ lies between $\bhat$ and $\b_0$. Similar to
approximating the order of $\hat{I}_3$ in the proof of Theorem
\ref{thm:ple}, the last line of the above equation is less than or
equal to
\begin{equation*}
n^{-2}O_p(n^2p_n^3)\|\bhat - \b_0\|^2 =
n^{-2}O_P(n^2p_n^3)O_P(p_n/n) = o_P(1)
\end{equation*}
by the conclusion of Theorem \ref{thm:ple}.
\thickspace\thickspace$\square$



\begin{lemma}\label{lemma:8}
Assuming Regularity Conditions (A) - (G), we have for each
$\b\in\Omega_n$,
$$ n^{-1}\|\nabla^2\hat{Q}_n(\b) - \nabla^2Q_n(\b)\| = O_P(p_n(h^{p+1}+c_n\lh)). $$
\end{lemma}

\noindent \emph{Proof of Lemma \ref{lemma:8}}. By Taylor's
expansion and Lemma \ref{lemma:order1},
\begin{equation*}
\begin{split}
n^{-1}&\frac{\partial}{\partial\beta_{nk}}(\nabla\hat{Q}_n(\b) - \nabla Q_n(\b))\\
&= n^{-1}  \biggl \{  \frac{\partial^3 Q_n(\b)} {\partial
\beta_{nk}\partial \b \partial \a_\b^T}(\ahat - \a_\b) +
\frac{\partial^2 Q_n(\b)}{\partial \b \partial \a_\b^T}\biggl (
\frac{\partial \ahat}{\partial \beta_{nk}} -
\frac{\partial \a_\b}{\partial \beta_{nk}} \biggr ) \\
&\;\;\;+ \biggl ( \frac{\partial \ahatp}{\partial \beta_{nk}} -
\frac{\partial \a_\b^\prime}{\partial \beta_{nk}} \biggr )
\frac{\partial Q_n(\b)}{\partial \a_\b} + (\ahatp - \a_\b^\prime)
\frac{\partial^2 Q_n(\b)}{\partial \a_\b \partial \beta_{nk}}
\biggr \}(1+o_P(1))
\end{split}
\end{equation*}
Hence, using Regularity Condition (C),
\begin{equation*}
\begin{split}
\bigg\|n^{-1}&\frac{\partial}{\partial\beta_{nk}}(\nabla\hat{Q}_n(\b) - \nabla Q_n(\b))\bigg\|\\
&= O(1)\cdot\bigg(\sup_{i}\|\ahat(U_i) - \a_\b(U_i)\| +
\sup_{i}\bigg\|\frac{\partial\ahat(U_i)}{\partial\beta_{nk}}- \frac{\partial\a_\b(U_i)}{\partial\beta_{nk}}\bigg\|\\
&+ \sup_{i}\|\ahatp(U_i) - \ahatp(U_i)\| +
\sup_{i}\bigg\|\frac{\partial\ahatp(U_i)}{\partial\beta_{nk}}
 - \frac{\a_\b^{\prime}(U_i)}{\partial\beta_{nk}}\bigg\|\bigg)\\
&= O_P(\sqrt{p_n}(h^{p+1} + c_n\lh)),
\end{split}
\end{equation*}
where the last line follows from Lemma \ref{lemma:order1}. Hence
$$n^{-1}\|\nabla^2\hat{Q}_n(\b) - \nabla^2Q_n(\b)\| =
O_P(p_n(h^{p+1} + c_n\lh)). \thickspace\thickspace\square$$

\begin{lemma}\label{lemma:1}
Under Regularity Conditions (A) - (G) and $p_n^4/n = o(1)$,
\begin{eqnarray*}
\|n^{-1}\nabla^2Q_n(\b_0)+I_n(\b_0)\| &=& o_P(p_n^{-1}), \\
\|n^{-1}\nabla^2\hat{Q}_n(\b_0)+I_n(\b_0)\| &=& o_P(p_n^{-1}) +
O_P(p_n(h^{p+1} + c_n\lh )).
\end{eqnarray*}
\end{lemma}

\noindent \emph{Proof of Lemma \ref{lemma:1}}. The first
conclusion follows from
\begin{equation*}
\begin{split}
&E_0 p_n^2 \|n^{-1} \nabla^2Q_n(\b_0)+I_n(\b_0)\|^2 \\
&=p_n^2n^{-2} E_0 \sum_{i,j=1}^{p_n}\biggl\{
\frac{\partial^2Q_n(\b_0)}{\partial\beta_{ni}\partial\beta_{nj}} -
E_0\frac{\partial^2Q_n(\b_0)}{\partial\beta_{ni}\partial\beta_{nj}}\biggr\}^2
 = O(p_n^4/n ) = o(1).
\end{split}
\end{equation*}
From this, triangle inequality immediately gives
$$ \|n^{-1}\nabla^2\hat{Q}_n(\b_0)+I_n(\b_0)\| =
o_P(p_n^{-1}) + \| n^{-1}\nabla^2(\hat{Q}_n(\b_0) - Q_n(\b_0))
\|.$$

The second equation then follows from Lemma \ref{lemma:8}.
$\thickspace\thickspace\square$



\begin{lemma}\label{lemma:4}
Assuming the conditions in Theorem \ref{thm:hyptest} and under the
null hypothesis $H_0$ as in the theorem,
$$ B_n^T(\ghatvec - \gvec_0) =  \frac{1}{n}B_n^T\{B_nI_n(\b_0)B_n^T\}^{-1}B_n^T\nabla Q_n(\b_0)
 + o_P(n^{-1/2}).$$
\end{lemma}

\noindent \emph{Proof of Lemma \ref{lemma:4}}. Since $B_nB_n^T =
I_{p_n-l}$, for each $\mathbf{v}\in \mathbb{R}^{p_n-l}$, we have
\begin{equation}\label{eqn:lemma4eqn1}
\|B_n^T\mathbf{v}\| \leq \|\mathbf{v}\|.
\end{equation}

Following the proof of Theorem \ref{thm:ple}, we have
$\|B_n^T(\ghatvec - \gvec)\| = O_P(\sqrt{p_n/n})$.
Following the proof of Theorem \ref{thm:asymnorm} and by
\mbox{Lemma \ref{lemma:order4}},
$$I_n(\b_0)B_n^T(\ghatvec - \gvec_0) =
n^{-1}\nabla Q_n(\b_0) + o_P( n^{-1/2} ).$$ Left-multiplying with
$B_n$ and using \mbox{equation (\ref{eqn:lemma4eqn1})}, the right
hand side of the above equation becomes $n^{-1}B_n\nabla Q_n(\b_0)
+ o_P(n^{-1/2})$. Hence,
$$B_n^T(\ghatvec - \gvec_0) = n^{-1}B_n^T(B_nI_n(\b_0)B_n^T)^{-1}B_n\nabla Q_n(\b_0)
+ o_P( n^{-1/2} ), $$ since $B_nI_n(\b_0)B_n^T$ has eigenvalues
uniformly bounded away from 0 and infinity, like $I_n(\b_0)$ does.
$\thickspace \square$

\end{document}